\def\titlep{Permutative representations of the Cuntz-Krieger algebras}

\documentclass[11pt]{article}
\usepackage{graphicx,ifthen}

 at12pt

\newcommand{\qed}{\hbox{\rule[-2pt]{3pt}{6pt}}}
\newcommand{\qedh}{\hfill\qed \\}
\newcommand{\vv}{\vspace{.3in}}

\newtheorem{Thm}{Theorem}[section]
\newtheorem{defi}[Thm]{Definition}
\newtheorem{lem}[Thm]{Lemma}

\newtheorem{prop}[Thm]{Proposition}
\newtheorem{cor}[Thm]{Corollary}

\newcommand{\ww}{\vv\noindent}

\def\cal#1{\mathcal #1}
\def\con{{\cal O}_{N}}
\def\edot{=1,\ldots,N}
\def\pr{{\it Proof.}\,\,}

\def\evp{eventually periodic}
\def\mnz{M_{N}(\{0,1\})}
\def\rnl{RN_{loc}(X)}
\def\coa{{\cal O}_{A}}
\def\co#1{{\cal O}_{#1}}
\def\ltn{l_{2}({\bf N})}

\def\mattwo#1#2#3#4{
\left(
\begin{array}{cc}
#1&#2\\
#3&#4\\
\end{array}
\right)}

\def\matthr#1#2#3#4#5#6#7#8#9{
\left(
\begin{array}{ccc}
#1&#2&#3\\
#4&#5&#6\\
#7&#8&#9\\
\end{array}
\right)}

\def\bfs{branching function system}
\def\bfset#1#2{{\rm BFS}_{#1}({#2})}

\addtocounter{footnote}{1}
\def\sftt#1{
\setcounter{equation}{0}
\addtocounter{footnote}{1}
\section{#1}
}

\def\ssft#1{\subsection{#1}}

\def\nset#1{\{1,\ldots,N\}^{#1}}

\def\brl{branching law}

\def\cka{Cuntz-Krieger algebra}

\def\mqb{\{(M_{i},q_{i},B_{i})\}_{i=1}^{N}}

\renewcommand{\S}{{\rm Sect.\,}}

\begin{document}
\title{{\titlep}}
%\titlerunning{Representations of the Cuntz-Krieger algebras}

\def\autherp{Katsunori  Kawamura}
\def\emailp{e-mail:kawamura@kurims.kyoto-u.ac.jp}

\def\addressp{College of Science and Engineering Ritsumeikan University,\\
1-1-1 Noji Higashi, Kusatsu, Shiga 525-8577, Japan.
}

%
%%%%%%%% Cut from here%%%%%%%%%%%%
%\input comm.txt
%
%\pagestyle{plain}
%\setcounter{page}{1}
%\setcounter{section}{0}
%%%%%%%% Cut end %%%%%%%%%%%%%%%

\author{\autherp}
%\institute{\addressp\\ \email{\emailp}}
%\thanks{{\it \addressp}}
%\authorrunning{K.Kawamura}

%\date{Received: *** / Accepted: ***}
%\communicated{****}

%
% Title
%
\begin{center}
{\Large \titlep}

\ww
\autherp
\footnote{\emailp}

\noindent
{\it \addressp}
%\quad \\
%\quad \\
\quad \\
%\fbox{{\Large \today}}

\end{center}

%\maketitle

%%%%%%%%%%%%%%%%%%%%%%%%%%%%%%%%%%%%
%
% Abstract
%
\begin{abstract}
We generalize permutative representations of the Cuntz algebras 
for the \cka\ $\coa$ for any $A$. 
We characterize cyclic permutative representations by notions of 
cycle and chain, and show their existence and uniqueness.
We show necessary and sufficient conditions for
their irreducibility and equivalence.
In consequence, we have a complete classification
of permutative representations of $\coa$.
Furthermore we show decomposition formulae and 
that the uniqueness of
irreducible decomposition holds for permutative
representation.
\end{abstract}

%
%%%%%%%%%%%%%%%%%%%%%%%%%%%%%%%%%%%
%
% Section 1
%
\sftt{Introduction}
\label{section:first}
Representation theory of the Cuntz-Krieger algebras
has an application to the Perron-Frobenius operator
of dynamical system \cite{PFO}.
For analysis of such representations,
it is necessary to research their fundamental property. 
On the other hand, permutative representations of the Cuntz algebras 
are completely classified by \cite{BJ,DaPi2}.
They are applied to representations of the fermion algebra \cite{AK1},
classification of theories of quantum string fields \cite{AK3}.
Further the classification of sectors are obtained by \brl s of 
permutative representations \cite{PE01}
and boson-fermion correspondence is explained by the branching law of 
permutative representation \cite{IWF01}.
In this article, we generalize permutative representations to the \cka\ 
$\coa$ for any $A$ and show their property, 
that is, existence, canonical basis, uniqueness, 
irreducibility, equivalence and irreducible decomposition.
The remarkable point is that {\it the uniqueness of irreducible decomposition holds} 
for permutative representations of $\coa$ for any $A$.
Therefore their decomposition formulae make sense.

Let $N\geq 2$ and $A$ be an $N\times N$ matrix
which has entries in $\{0,1\}$ and has no rows or columns 
identically equal to zero, and
let $s_{1},\ldots,s_{N}$ be canonical generators of $\coa$.
%
% Definition 1.1
%
\begin{defi}
\label{defi:chara}
Let $({\cal H},\pi)$ be a representation of $\coa$.
\begin{enumerate}
%(i)
\item
$({\cal H},\pi)$ is a {\bf permutative representation} of $\coa$ with phases
if there is a complete orthonormal basis
$\{e_{n}\}_{n\in\Lambda}$ of ${\cal H}$,
a family $\{\Lambda_{i}\}_{i=1}^{N}$ of subsets of $\Lambda$
and $(z_{i,n},m_{i,n})\in U(1)\times \Lambda$ such that
%
% Equation 1.1
%
\begin{equation}
\label{eqn:first}
\pi(s_{i})e_{n}=z_{i,n}e_{m_{i,n}}\quad (n\in\Lambda_{i}),
\quad \pi(s_{i})e_{n}=0\quad(\mbox{otherwise}).
\end{equation}
Especially, we call $({\cal H},\pi)$ by a permutative representation 
of $\coa$ when $z_{i,n}=1$ for any $(i,n)$.
%(i)
\item
$({\cal H},\pi)$ is $P(J;z)$ for $J=(j_{l})_{l=1}^{p}\in
\nset{p}$ and $z\in U(1)$ 
if there is a cyclic unit vector $\Omega\in{\cal H}$ such that
$\pi(s_{j_{1}}\cdots s_{j_{p}})\Omega=z\Omega$ and
$\{\pi(s_{j_{l}}\cdots s_{j_{p}})\Omega\}_{l=1}^{p}$
is an orthonormal family.
Especially, we denote $P(J)\equiv P(J;1)$.
%(ii)
\item
$({\cal H},\pi)$ is $P(J)$ for $J=(j_{n})_{n\in{\bf N}}\in\nset{\infty}$
if there is a cyclic unit vector $\Omega\in{\cal H}$
such that $\{\pi(s_{j_{n}}^{*}\cdots s_{j_{1}}^{*})\Omega\}_{n\in {\bf N}}$
is an orthonormal family.
\end{enumerate}
$\Omega$ in both (ii) and (iii) is called the GP vector
of $({\cal H},\pi)$.
\end{defi}
%
% Theorem 1.2
%
\begin{Thm}
\label{Thm:matwo}
If $({\cal H},\pi)$ is a permutative representation of $\coa$ with phases,
then the following holds:
\begin{enumerate}
% (i)
\item
$({\cal H},\pi)$ is decomposed into the direct sum of cyclic permutative 
representations with phases uniquely up to unitary equivalence.
%
% (ii)
\item
If $({\cal H},\pi)$ is cyclic, then
$({\cal H},\pi)$ is either the case (i) or (ii) in Definition \ref{defi:chara}.
%
% (iii)
\item
$P((j_{n})_{n=1}^{p};c)$ {\rm ({\it resp}. $P((k_{n})_{n\in{\bf N}})$)}
is irreducible if and only if there is no $\sigma\in {\bf Z}_{p}\setminus
\{id\}$ such that $(j_{\sigma(1)},\ldots,j_{\sigma(p)})
=(j_{1},\ldots,j_{p})$ {\rm (}{\it resp}.
there is no $(q,n_{0}) \in {\bf N}\times {\bf N}$
such that $k_{n+q}=k_{n}$ for each $n\geq n_{0}$.{\rm )}
%
%(iv)
\item
$P((j_{n})_{n=1}^{p};c)\not \sim P((k_{n})_{n\in{\bf N}})$
for any $((j_{n})_{n=1}^{p},c)$ and  $(k_{n})_{n\in{\bf N}}$.
%(v)
\item
$P((j_{n})_{n=1}^{p};c)\sim P((j_{n}^{'})_{n=1}^{p^{'}};c^{'})$
if and only if $p=p^{'}$, $c=c^{'}$
and there is $\sigma\in{\bf Z}_{p}$ such that
$j_{\sigma(n)}^{'}=j_{n}$ for each $n=1,\ldots,p$.
%(vi)
\item
$P((k_{n})_{n\in{\bf N}}) \sim P((k_{n}^{'})_{n\in{\bf N}})$
if and only if there is $(q,n_{0})\in {\bf Z}\times
{\bf N}$ such that $k_{n+q}=k^{'}_{n}$ for each $n\geq n_{0}$.
\end{enumerate}
\end{Thm}

%
%  2
In $\S$\ref{section:second}, we prepare multiindices
associated with a matrix $A$ and
introduce $A$-\bfs s and show their properties.
%
% 3
In $\S$\ref{section:third}, 
the existence of cyclic representations appearing in Theorem \ref{Thm:matwo} (ii)
is shown for each multiindex in $\S$\ref{section:second}.
We show the construction of the canonical basis
of a given permutative representation. 
By using this, we show uniqueness.
%
% 4
In $\S$\ref{section:fourth}, we show 
irreducibility and equivalence of them.
%
% 5
In $\S$\ref{section:fifth}, we show decomposition 
formulae of permutative representations. 
Theorem \ref{Thm:matwo} is shown here.
%
% 6
In $\S$\ref{section:sixth}, we show states and
spectrums of $\coa$ associated with
permutative representations.
%
% 7
In $\S$\ref{section:seventh}, we show examples.

%%%%%%%%%%%%%%%%%%%%%%%%%%%%%%%%%%%%%%%%%%%
%
% Section 2
%
\sftt{$A$-\bfs s}
\label{section:second}
%
%%%%%%%%%%%%%%%%%%%%%%%%%%%%%%%%%%%%
%
% subsection 2.1
%
\ssft{Multiindices}
\label{subsection:secondone}
Let $\nset{0}\equiv \{0\}$, 
$\nset{k}\equiv\{(j_{l})_{l=1}^{k}:j_{l}\edot,\,l=1,\ldots,k\}$
for $k\geq 1$ and
$\nset{\infty}\equiv\{(j_{n})_{n\in {\bf N}}:j_{n}\in \nset{},\,
n\in{\bf N}\}$, $\nset{*}$ $\equiv \coprod_{k\geq 0}\nset{k}$,
$\nset{*}_{1}\equiv\coprod_{k\geq 1}\nset{k}$,
$\nset{\#}\equiv$ $\nset{*}_{1}$ $\sqcup \nset{\infty}$.
For $J\in \nset{\#}$, the {\it length} $|J|$ of $J$ is defined by
$|J|\equiv k$ when $J\in \nset{k}$.
For $J_{1},J_{2}\in\nset{*}$ and $J_{3}\in\nset{\infty}$,
$J_{1}\cup J_{2}\equiv(j_{1},\ldots,j_{k},j_{1}^{'},\ldots,j_{l}^{'})$,
$J_{1}\cup J_{3}\equiv(j_{1},\ldots,j_{k},j_{1}^{''},j_{2}^{''},\ldots)$
when $J_{1}=(j_{1},\ldots,j_{k})$, $J_{2}=(j_{1}^{'},\ldots,j_{l}^{'})$
and $J_{3}=(j_{n}^{''})_{n\in {\bf N}}$.
Especially, we define $J\cup \{0\}=\{0\}\cup J=J$
for $J\in \nset{\#}$ and $(i,J)\equiv (i)\cup J$ for convenience.
For $J\in\nset{*}$, $J^{k}\equiv J\cup\cdots\cup J$ ($k$ times) 
and $J^{\infty}\equiv J\cup J\cup J\cup\cdots\in \nset{\infty}$.
For $J=(j_{1},\ldots,j_{k})$ and $\tau\in {\bf Z}_{k}$,
define $\tau(J)\equiv (j_{\tau(1)},\ldots,j_{\tau(k)})$.

For $N\geq 2$, let $\mnz$ be the set of all $N\times N$ matrices in
which have entries in $\{0,1\}$ and has no rows or columns identically 
equal to zero.
$A=(a_{ij})$ is {\it full} if $a_{ij}=1$ for each $i,j\edot$.
For $A=(a_{ij})\in\mnz$, define
\[\nset{*}_{A}\equiv \coprod_{k\geq 0}\nset{k}_{A},\quad
\nset{*}_{A,c}\equiv \coprod_{k\geq 1}\nset{k}_{A,c},\]
\[\nset{\infty}_{A}\equiv\{(j_{n})_{n\in{\bf N}}\in
\nset{\infty}:a_{j_{n-1}j_{n}}= 1,n\geq 2\},\]
\[\nset{\#}_{A,c}\equiv\nset{*}_{A,c}\sqcup\nset{\infty}_{A}\]
where $\nset{0}_{A}\equiv \{0\}$, $\nset{1}_{A}\equiv \nset{}$,
$\nset{k}_{A}\equiv\{(j_{i})_{i=1}^{k}\in\nset{k}:a_{j_{i-1}j_{i}}=1,
\,i=2,\ldots,k \}$, 
$\nset{k}_{A,c}\equiv\{(j_{i})_{i=1}^{k}\in\{1,\ldots,N\}^{k}_{A}:
a_{j_{k}j_{1}}=1 \}$ for $k\geq 2$.

$J\in\nset{*}_{1}$
is {\it periodic} if
there are $m\geq 2$
and $J_{0}\in \nset{*}_{1}$
such that
$J=J_{0}^{m}$.
For
$J_{1},J_{2}\in \nset{*}_{1}$,
$J_{1}\sim J_{2}$
if 
there are $k\geq 1$
and $\tau\in {\bf Z}_{k}$
such that
$|J_{1}|=|J_{2}|=k$
and
$\tau(J_{1})=J_{2}$.
For
$(J,z),
(J^{'},z^{'})\in \nset{*}_{1}\times U(1)$,
$(J,z)\sim
(J^{'},z^{'})$
if 
$J\sim J^{'}$
and $z=z^{'}$
where
$U(1)\equiv \{z\in {\bf C}:|z|=1\}$.
$J\in\nset{\infty}$
is {\it \evp} if 
there
are
$J_{0},J_{1}\in\nset{*}_{1}$
such that
$J=J_{0}\cup J_{1}^{\infty}$.
Especially,
if
$J\in\nset{\infty}_{A}$,
then
$J_{0}\in\nset{*}_{A}$
and
$J_{1}\in\nset{*}_{A,c}$
in the above.
For
$J_{1},J_{2}\in\nset{\infty}$,
$J_{1}\sim J_{2}$
if there
are $J_{3},J_{4}\in\nset{*}$
and 
$J_{5}\in\nset{\infty}$
such that
$J_{1}=J_{3}\cup J_{5}$
and
$J_{2}=J_{4}\cup J_{5}$.
If
$J\in\nset{\infty}_{A}$ is \evp,
then
there is 
$J_{1}\in\nset{*}_{A,c}$
such that
$J\sim J_{1}^{\infty}$.
For
$J,J^{'}\in\nset{\#}$,
$J\sim J^{'}$
if 
$J,J\in
\nset{*}_{1}$
and 
$J\sim J^{'}$,
or
$J,J\in
\nset{\infty}$
and 
$J\sim J^{'}$.

For
$J_{1}=(j_{1},\ldots,j_{k}), J_{2}=(j_{1}^{'},\ldots,j_{k}^{'})\in\nset{k}$,
$J_{1}\prec J_{2}$
if $\sum_{l=1}^{k}(j_{l}^{'}-j_{l})N^{k-l}\geq 0$.
$J\in\nset{*}_{1}$ is {\it minimal} if $J\prec J^{'}$ for each 
$J^{'}\in\nset{*}_{1}$ such that $J\sim J^{'}$.
Especially, any element in $\nset{}$
is non periodic and minimal.
Define
\[<1,\ldots,N>^{*}_{A}\equiv \{J\in\nset{*}_{A,c}:J\mbox{ is minimal}\},\]
\[<1,\ldots,N>^{\infty}_{A}\equiv \nset{\infty}_{A}/\!\!\sim,\]
\[[1,\ldots,N]^{*}_{A}\equiv\{J\in<1,\ldots,N>^{*}_{A}:
J\mbox{ is non periodic}\},\]
\[[1,\ldots,N]^{\infty}_{A}\equiv\{[J]\in<1,\ldots,N>^{\infty}_{A}:
J\mbox{ is non \evp}\}\]
where
$[J]\equiv\{J^{'}\in\nset{\infty}_{A}:J\sim J^{'}\}$.
Then $[1,\ldots,N]^{*}_{A}$ is in one-to-one correspondence
with the set of all equivalence classes of 
non periodic elements in $\nset{*}_{A,c}$.
We denote an element $[K]$ of both $<1,\ldots,N>^{\infty}_{A}$
and $[1,\ldots,N]^{\infty}_{A}$ by a representative element $K$
if there is no ambiguity.
Define
%
% Equation 2.1
%
\begin{equation}
\label{eqn:alldef}
\left\{
\begin{array}{l}
<1,\ldots,N>^{\#}_{A}\equiv 
<1,\ldots,N>^{*}_{A}\sqcup
<1,\ldots,N>^{\infty}_{A},\\
\\
\mbox{$[1,\ldots,N\,]^{\#}_{A}$}\equiv 
[1,\ldots,N]^{*}_{A}\sqcup
[1,\ldots,N]^{\infty}_{A}.
\end{array}
\right.
\end{equation}
%

%%%%%%%%%%%%%%%%%%%%%%%%%%%%%%%%%%%%%%%%%%
%
% subsection 2.2
%
\ssft{$A$-\bfs s}
\label{subsection:secondtwo}
We denote the set of injective maps from $X$ to $Y$ by $RN(X,Y)$ and
define $RN_{loc}(X,Y)\equiv\bigcup_{X_{0}\subset X}RN(X_{0},Y)$.
We simply denote $RN(X)\equiv RN(X,X)$.
For $f\in \rnl$, we denote the domain and the range of 
$f$ by $D(f)$ and $R(f)$, respectively.
$\rnl$ and $RN(X)$ is a groupoid and a semigroup 
by composition of maps, respectively.
We denote $X\times Y$ and $X\cup Y$, the direct product 
and the direct sum of $X$ and $Y$ as sets, respectively.
For $f\in RN(X_{1},Y_{1})$ and $g\in RN(X_{2},Y_{2})$,
$f\oplus g\in RN(X_{1}\cup X_{2},Y_{1}\cup Y_{2})$
is defined by $(f\oplus g)|_{X_{1}}\equiv f$,
$(f\oplus g)|_{X_{2}}\equiv g$.

%
% Definition 2.1
% 
\begin{defi}
\label{defi:abfs}
For $A=(a_{ij})\in \mnz$, $f=\{f_{i}\}_{i=1}^{N}$
is an $A$-\bfs\ on a set $X$ if $f$ satisfies the following:
\begin{enumerate}
%(i)
\item
There is a family $\{D(f_{i})\}_{i=1}^{N}$ of subsets of $X$
such that $f_{i}$ is an injective map from $D(f_{i})$ to $X$
with the image $R(f_{i})$ for each $i\edot$.
%(ii)
\item
$R(f_{i})\cap R(f_{j})=\emptyset$ when $i\ne j$.
%(iii)
\item
$D(f_{i})=\coprod_{j:a_{ij}=1}R(f_{j})$ for each $i\edot$.
%(iv)
\item
$X= \coprod_{i=1}^{N}R(f_{i})$.
\end{enumerate}
Especially, if $A$ is full, then we call $A$-\bfs\ by ($N$-) \bfs\ simply.
We denote the set of all $A$-\bfs s, \bfs s on $X$ 
by $\bfset{A}{X}$, $\bfset{N}{X}$, respectively.
\end{defi}

\noindent
By definition, $\bfset{A}{X}\ne \emptyset$ if and only if $\#X=\infty$.
Hence we assume $\#X=\infty$ in this article.
$F$ is the {\it coding map} of $f=\{f_{i}\}_{i=1}^{N}\in\bfset{A}{X}$
if $F$ is a map on $X$ such that $(F\circ f_{i})(x)=x$
for each $x\in X$ and $i\edot$.
For $f=\{f_{i}\}_{i=1}^{N}\in\bfset{A}{X}$ and
$g=\{g_{i}\}_{i=1}^{N}\in\bfset{A}{Y}$, 
$f\sim g$ if there is a bijection $\varphi$ from $X$ to $Y$
such that $\varphi\circ f_{i}\circ \varphi^{-1}=g_{i}$ for $i\edot$.
For $f=\{f_{i}\}_{i=1}^{N}\in\bfset{A}{X}$ and
$g=\{g_{i}\}_{i=1}^{N}\in\bfset{A}{Y}$,
we denote $f\oplus g\equiv \{f_{i}\oplus g_{i}\}_{i=1}^{N}
\in\bfset{A}{X\cup Y}$.
Let $\{X_{\omega}\}_{\omega\in\Xi}$ be a family of sets.
For $f^{[\omega]}=\{f_{i}^{[\omega]}\}_{i=1}^{N}\in\bfset{A}{X_{\omega}}$
for $\omega\in\Xi$,
$f$ is the {\it direct sum} of $\{f^{[\omega]}\}_{\omega\in\Xi}$
if $f=\{f_{i}\}_{i=1}^{N}\in\bfset{A}{X}$ for
a set $X\equiv \coprod_{\omega\in\Xi}X_{\omega}$ which
is defined by $f_{i}(n)\equiv f^{[\omega]}_{i}(n)$
when $n\in X_{\omega}$ for $i\edot$ and $\omega\in\Xi$.
For $f\in \bfset{A}{X}$, $f=\bigoplus_{\omega\in\Xi}
f^{[\omega]}$ is a {\it decomposition} of $f$
into a family $\{f^{[\omega]}\}_{\omega\in\Xi}$
if there is a family $\{X_{\omega}\}_{\omega\in\Xi}$
of subsets of $X$ such that $f$ is the direct sum of $\{
f^{[\omega]}\}_{\omega\in\Xi}$.

For $f=\{f_{i}\}_{i=1}^{N}\in\bfset{A}{X}$,
define $f_{J}\equiv f_{j_{1}}\circ\cdots\circ f_{j_{k}}$
when $J=(j_{l})_{l=1}^{k}$ and define $f_{0}\equiv id$.
When we denote $f_{i}(x)$, we assume $x\in D(f_{i})$ automatically.
Define ${\cal C}_{x}\equiv\{f_{J}(x)\in X:J\in\nset{*}_{A}\, 
s.t.\, x\in D(f_{J})\}\cup\{F^{n}(x)\in X:n\in{\bf N}\}$
where $F$ is the coding map of $f$.
For $x\in X$, if $y\in {\cal C}_{x}$, then ${\cal C}_{y}={\cal C}_{x}$.
For $f\in \bfset{A}{X}$, $f$ is {\it cyclic} if there is $x\in X$ 
such that ${\cal C}_{x}=X$.
$\{x_{i}\}_{i=1}^{k}\subset X$ is a {\it cycle} of $f$ by $J=(j_{l})_{l=1}^{k}$
if $f_{j_{l}}(x_{l+1})=x_{l}$ for $l=1,\ldots,k-1$ and $f_{j_{k}}(x_{1})=x_{k}$.
$\{x_{n}\}_{n\in{\bf N}}\subset X$ is a {\it chain} of $f$ by 
$J=(j_{n})_{n\in{\bf N}}$ if $f_{j_{n-1}}(x_{n})=x_{n-1}$ for each $n\geq 2$.
For each $x\in X$ and $f\in \bfset{A}{X}$,
$f|_{{\cal C}_{x}}\in \bfset{A}{{\cal C}_{x}}$
and $f|_{{\cal C}_{x}}$ is cyclic.
For $A\in\mnz$, let $f\in \bfset{A}{X}$.
If $f$ is cyclic, then $f$ has either only a cycle or a chain.
If $f$ is cyclic and has two chains $\{x_{n}\}_{n\in {\bf N}}$
and $\{y_{n}\}_{n\in {\bf N}}$, then
there are $p$ and $M\geq 0$ such that
$x_{p+n}=y_{n}$ or $x_{n}=y_{n+p}$ for each $n> M$.
For any $f\in\bfset{A}{X}$, there is a
decomposition $X=\coprod_{\lambda\in\Lambda}
X_{\lambda}$ such that $f|_{X_{\lambda}}$
is cyclic for each $\lambda\in \Lambda$.

%
% Definition 2.3
%
\begin{defi}
\label{defi:pjf}
For $A\in\mnz$, let $f\in \bfset{A}{X}$.
\begin{enumerate}
%(i)
\item
For $J\in\nset{*}_{A,c}$ {\rm ({\it resp}. $J\in\nset{\infty}_{A}$)},
$f$ has a $P(J)$-component if $f$ has a cycle ({\it resp}. a chain) by $J$.
%(ii)
\item
For $J\in\nset{\#}_{A,c}$, $f$ is $P(J)$
if $f$ is cyclic and has a $P(J)$-component.
\end{enumerate}
\end{defi}

\noindent
For $J,J^{'}\in\nset{\#}_{A,c}$, assume that
$f$ and $f^{'}$ are $P(J)$ and $P(J^{'})$, respectively.
Then $f\sim f^{'}$ if and only if $J\sim J^{'}$.
This follows from the uniqueness of cycle and chain up to
equivalence.
In consequence, the following holds:
%
% Lemma 2.4
%
\begin{lem}
\label{lem:decothm}
For any $A\in\mnz$ and $f\in\bfset{A}{X}$,
there is a decomposition $X=\coprod_{\lambda\in\Lambda}X_{\lambda}$
where $f|_{X_{\lambda}}$ is $P(J_{\lambda})$
for some $J_{\lambda}\in\nset{\#}_{A,c}$ for each $\lambda\in\Lambda$.
This decomposition is unique up to equivalence of \bfs s. 
\end{lem}

\noindent
We can simply describe the statement in Lemma \ref{lem:decothm}
as follows:
\[f\sim\bigoplus_{J\in<1,\ldots,N>^{\#}_{A}}P(J)^{\oplus \nu_{J}}\]
where $\nu_{J}$ is the multiplicity of $P(J)$ in $f$.

%%%%%%%%%%%%%%%%%%%%%%%%%%%%%%%%%%%%%%%%%%%%%
%
% Subsection 2.3
%
\ssft{Construction of $A$-\bfs}
\label{subsection:secondthree}
Fix $A=(a_{ij})\in\mnz$.
Define
%
% Equation  2.2
%
\begin{equation}
\label{eqn:tree}
{\cal T}(A;j)\equiv \coprod_{k\geq 1}{\cal T}^{(k)}(A;j),\quad
{\cal T}(j;A)\equiv \coprod_{k\geq 1}{\cal T}^{(k)}(j;A)
\end{equation}
where
${\cal T}^{(k)}(A;j)\equiv
\{(j_{1},\ldots,j_{k})\in\nset{k}_{A}:a_{j_{k}j}=1\}$
and
${\cal T}^{(k)}(j;A)\equiv \{(j_{1},\ldots,j_{k})\in\nset{k}_{A}:
a_{jj_{1}}=1\}$.
For $J=(j_{1},\ldots,j_{k})$, 
define $J_{l}\equiv (j_{l},\ldots,j_{k})$ for $l=1,\ldots,k$ and 
%
% Equation 2.3 
%
\begin{equation}
\label{eqn:lambda}
\Lambda(A,J)\equiv\Lambda_{1}(A,J)\sqcup\Lambda_{2}(A,J)\sqcup\Lambda_{3}(A,J),
\end{equation}
$\Lambda_{1}(A,J)\equiv\{J_{l}:1\leq l\leq k\}$,
$\Lambda_{2}(A,J)\equiv\Lambda_{2,1}(A,J)\sqcup\cdots\sqcup\Lambda_{2,k}(A,J)$,
$\Lambda_{2,l}(A,J)\equiv\{(j,J_{l}):j\in {\cal T}^{(1)}(A;j_{l}),\,
j\ne j_{l-1}\}$ for $l=1,\ldots,k$,
$\Lambda_{3}(A,J)\equiv\coprod_{(j,J_{l})\in\Lambda_{2}(A,J)}
{\cal T}(A;j)\times \{(j,J_{l})\}$
where $j_{0}\equiv j_{k}$.
Then the following holds:
%
% Lemma 2.4
%
\begin{lem}
\label{lem:existcycle}
Let $\{D(f_{i})\}_{i=1}^{N}$ be a family
of subsets of $\Lambda(A,J)$ by
\[D(f_{i})\equiv{\cal T}(i;A)\cap \Lambda(A,J)\quad(i\edot)\]
and $f=\{f_{i}\}_{i=1}^{N}$ be a family of maps by
$f_{i}:D(f_{i})\to \Lambda(A,J)$;
\[f_{i}(J^{'})\equiv J_{k}\quad (J^{'}= J \mbox{ and }i=j_{k}),
\quad f_{i}(J^{'})\equiv(i,J^{'})\quad(\mbox{otherwise})
\]
for $i\edot$.
Then $f$ is an $A$-\bfs\ on $\Lambda(A,J)$ and $f$ is $P(J)$.
\end{lem}

For $J=(j_{n})_{n\in {\bf N}}\in \nset{\infty}_{A}$, define
%
% Equation 2.4
% 
\begin{equation}
\label{eqn:lambdacha}
\Lambda(A,J)\equiv
\Lambda_{1}(A,J)\sqcup
\Lambda_{2}(A,J)\sqcup
\Lambda_{3}(A,J)
\end{equation}
where $\Lambda_{1}(A,J)\equiv{\bf N}$,
$\Lambda_{2}(A,J)\equiv\coprod_{m\in{\bf N}}\Lambda_{2,m}(A,J)$,
$\Lambda_{2,1}(A,J)\equiv\{(j,1):j\in {\cal T}^{(1)}(A;j_{1})\}$,
$\Lambda_{2,m}(A,J)\equiv\{(j,m):j\in {\cal T}^{(1)}(A;j_{m}),\,j\ne 
j_{m-1}\}$ for $m\geq 2$ and 
$\Lambda_{3}(A,J)\equiv\coprod_{(j,m)\in
\Lambda_{2}(A,J)}{\cal T}(A;j)\times \{(j,m)\}$.
Then the following holds:
%
% Lemma 2.5
%
\begin{lem}
\label{lem:existchain}
Let $\{D(f_{i})\}_{i=1}^{N}$ be a family of subsets of 
$\Lambda(A,J)$ by
\[D(f_{i})\equiv\{m\in{\bf N}: a_{ij_{m}}=1\}\sqcup
({\cal T}(i;A)\times {\bf N})\cap \Lambda(A,J)\]
and let $f=\{f_{i}\}_{i=1}^{N}$ be a family of maps by
$f_{i}:D(f_{i})\to \Lambda(A,J)$;
\[
\left\{
\begin{array}{l}
f_{i}(m)\equiv 
\left\{
\begin{array}{ll}
m-1\quad&(i=j_{m-1}\mbox{ and }m\geq 2),\\
&\\
(i,m)\quad&(\mbox{otherwise})\\
\end{array}
\right.\quad
(m\in\Lambda_{1}(A,J)\cap D(f_{i})),\\
\\
f_{i}(J^{'},m)\equiv(\{i\}\cup J^{'},\, m)\quad\quad\quad
((J^{'},m)\in(\Lambda_{2}(A,J)\sqcup\Lambda_{3}(A,J))\cap D(f_{i})).
\end{array}
\right.
\]
Then $f$ is an $A$-\bfs\ on $\Lambda(A,J)$ and $f$ is $P(J)$.
\end{lem}

\noindent
By Lemma \ref{lem:existcycle} and Lemma \ref{lem:existchain},
the following holds:
%
% Lemma 2.6
%
\begin{lem}
\label{lem:existbfs}
For each $A\in\mnz$ and $J\in\nset{\#}_{A,c}$,
there is an element in $\bfset{A}{{\bf N}}$ which is $P(J)$.
\end{lem}

%%%%%%%%%%%%%%%%%%%%%%%%%%%%%%%%%%%%%%%%%%%
%
% Section 3
%
\sftt{Existence and uniqueness}
\label{section:third}
For $A=(a_{ij})\in \mnz$, $\coa$ is {\it the \cka\ by $A$} 
if $\coa$ is a C$^{*}$-algebra
which is universally generated by
partial isometries $s_{1},\ldots,s_{N}$ satisfying
$s_{i}^{*}s_{i}=\sum_{j=1}^{N}a_{ij}s_{j}s_{j}^{*}$
for $i\edot$ and $\sum_{i=1}^{N}s_{i}s_{i}^{*}=I$ \cite{CK}.
Especially, $\coa$ is the Cuntz algebra $\con$ when $A$ is full.
For $g=(z_{1},\ldots,z_{N})\in T^{N}(\equiv U(1)^{N})$,
define $\alpha_{g}\in {\rm Aut}\coa$
by $\alpha_{g}(s_{i})\equiv z_{i}s_{i}$ for $i\edot$.
We denote the canonical $U(1)$-gauge action on $\coa$ by $\gamma$.
For $J=(j_{l})_{l=1}^{k}$, we denote $s_{J}\equiv s_{j_{1}}\cdots s_{j_{k}}$
and $s_{J}^{*}\equiv s_{j_{k}}^{*}\cdots s_{j_{1}}^{*}$.
When $J\in\nset{*}$, $s_{J}\ne 0$ if and only if $J\in\nset{*}_{A}$.
In this article, a representation always 
means a unital $^{*}$-representation on a complex Hilbert space.
$({\cal H}_{1},\pi_{1})\sim ({\cal H}_{2},\pi_{2})$
means the unitary equivalence between
two representations $({\cal H}_{1},\pi_{1})$ and
$({\cal H}_{2},\pi_{2})$ of $\coa$.

For $f=\{f_{i}\}_{i=1}^{N}\in\bfset{A}{X}$,
a representation $(l_{2}(X),\pi_{f})$ of $\coa$ is defined by
$\pi_{f}(s_{i})e_{n}=\chi_{D(f_{i})}(n)\cdot e_{f_{i}(n)}$
for $i\edot,\,n\in X$
where $\chi_{D(f_{i})}$ is the characteristic function on $D(f_{i})$.
Then $(l_{2}(X),\pi_{f})$ is a permutative representation.
Further any permutative representation of $\coa$ is equivalent to
$(l_{2}(X),\pi_{f})$ for a suitable $f$.
We denote $(l_{2}(X),\pi_{f})$ by $\pi_{f}$ simply.
For $f \in\bfset{A}{X}$ and $g\in\bfset{A}{Y}$,
$\pi_{f\oplus g}\sim\pi_{f}\oplus\pi_{g}$.
If $f\sim g$, then $\pi_{f}\sim\pi_{g}$.
If $f$ is cyclic, then $\pi_{f}$ is.
Recall $P(J)$ in Definition \ref{defi:abfs}.
If $f$ contains a $P(J)$-component for $J\in\nset{\#}_{A,c}$,
then $\pi_{f}$ does.
For $J\in\nset{\#}_{A,c}$, if $f$ is $P(J)$, then $\pi_{f}$ is $P(J)$.
For a representation $({\cal H},\pi)$ which is $P(J)$
and an automorphism $\varphi$, 
we denote $({\cal H},\pi\circ \varphi)$ by $P(J)\circ \varphi$.
By Definition \ref{defi:chara}, the following holds:

%
% Lemma 3.1
%
\begin{lem}
\label{lem:guageaction}
\begin{enumerate}
%(i)
\item
For $J=(j_{1},\ldots,j_{k})\in\nset{k}_{A,c}$,
$g=(z_{i})_{i=1}^{N}\in T^{N}$ and $w\in U(1)$,
$P(J;w)\circ \alpha_{g}\sim P(J;wz_{J})$
where $z_{J}\equiv z_{j_{1}}\cdots z_{j_{k}}$.
Especially, for $z,w\in U(1)$,
$P(J;w)\circ \gamma_{z}\sim P(J;wz^{k})$ and
$P(J)\circ \gamma_{z}\sim P(J;z^{k})$.
%(ii)
\item
For $J\in\nset{\infty}_{A}$ and $g\in T^{N}$,
$P(J)\circ \alpha_{g}\sim P(J)$.
\end{enumerate}
\end{lem}

\noindent
By this, Lemma \ref{lem:existbfs} and
Lemma \ref{lem:guageaction}, we have the following:
%
% Lemma 3.2
%
\begin{lem}
\label{lem:bfsandrep}
Let $A\in\mnz$.
\begin{enumerate}
%(i)
\item
For each $J\in\nset{\#}_{A,c}$,
there exists a representation of $\coa$ which is $P(J)$.
%(ii)
\item
For each $J\in\nset{*}_{A,c}$ and $z\in U(1)$,
there exists a representation of $\coa$ which is $P(J;z)$.
\end{enumerate}
\end{lem}

%
% Proposition 3.3
%
\begin{prop}
\label{prop:decopro}
For any permutative representation $({\cal H},\pi)$ of $\co{A}$,
there is a family $\{({\cal H}_{\lambda},\pi_{\lambda})\}_{\lambda\in\Lambda}$
of cyclic permutative representations of $\coa$
such that $({\cal H},\pi)=\bigoplus_{\lambda\in
\Lambda}({\cal H}_{\lambda},\pi_{\lambda})$.
Furthermore $({\cal H}_{\lambda},\pi_{\lambda})$
is $P(J_{\lambda})$ for some $J_{\lambda}\in\nset{\#}_{A,c}$
for each $\lambda\in\Lambda$.
\end{prop}
%
% Proof
%
\pr
By Lemma \ref{lem:decothm} and Lemma \ref{lem:bfsandrep}, it holds.
\qedh

\noindent
By definition of $P(J)$, Lemma \ref{lem:existcycle}
and Lemma \ref{lem:existchain}, we obtain the following:
%
% Lemma 3.4
%
\begin{lem}
\label{lem:cancy}
For $J\in\nset{\#}_{A}$, let $({\cal H},\pi)$ be $P(J)$
with the GP vector $\Omega$.
Then the following holds:
\begin{enumerate}
%(i)
\item
When $J\in\nset{*}_{A,c}$, define 
$e_{x}\equiv \pi(s_{x})\Omega$ for $x\in\Lambda(A,J)$
where $\Lambda(A,J)$ is in (\ref{eqn:lambda}).
Then $\{e_{x}\}_{x\in\Lambda(A,J)}$
is a complete orthonormal basis of ${\cal H}$.
%(ii)
\item
When $J=(j_{n})_{n\in {\bf N}}\in\nset{\infty}_{A}$, define
\[e_{n}\equiv \pi(s_{J[1,n]})^{*}\Omega\quad(n\in{\bf N}),\quad
e_{J^{'},n}\equiv\pi(s_{J^{'}}) e_{n}\quad(J^{'}\in {\cal T}(A;j_{n})
\setminus\{j_{n-1}\}),\]
\[\widetilde{\Lambda}(A,J)\equiv {\cal T}(A;j_{1})\sqcup 
\coprod_{n\geq 1}{\cal T}^{out}_{n}(A;J),\quad 
{\cal T}^{out}_{n}(A;J) \equiv {\cal T}(A;j_{n})\setminus\{j_{n-1}\}\]
where $J[1,n]\equiv (j_{1},\ldots,j_{n})$
and ${\cal T}(A;j_{n})$ is in (\ref{eqn:tree}).
Then $\{e_{x}\}_{x\in\widetilde{\Lambda}(A,J)}$
is a complete orthonormal basis of ${\cal H}$.
\end{enumerate}
\end{lem}
By Lemma \ref{lem:guageaction},  Lemma \ref{lem:bfsandrep} 
and Lemma \ref{lem:cancy},  the following holds.

%
% Proposition 3.5
%
\begin{prop}
\label{prop:uniqeness}
\begin{enumerate}
%(i)
\item
For $J\in\nset{\#}_{A,c}$, $P(J)$ exists uniquely
up to unitary equivalence.
%(ii)
\item
For $J\in\nset{*}_{A,c}$ and $z\in U(1)$, 
$P(J;z)$ exists uniquely up to unitary equivalence.
\end{enumerate}
\end{prop}
By Proposition \ref{prop:uniqeness},
both symbols $P(J)$ and $P(J;z)$ make sense
as an equivalence class of representations.

%%%%%%%%%%%%%%%%%%%%%%%%%%%%%%%%%%%%%%
%
% section 4
%
\sftt{Irreducibility and equivalence}
\label{section:fourth}
Let $A\in\mnz$.
%%%%%%%%%%%%%%%%%%%%%%%%%%%%%%%%%%%%%%
%
% subsection 4.1
%
\ssft{Sufficient condition
of irreducibility}
\label{subsection:fourthone}
%
% Lemma 4.1
%
\begin{lem}
\label{lem:lemirre}
Let $({\cal H},\pi)$ be $P(J)$ with the GP vector $\Omega$ for
$J=(j_{l})_{l=1}^{k}\in\nset{k}_{A,c}$ and
$\Omega_{l}\equiv \pi(s_{j_{l}}\cdots s_{j_{k}})\Omega$ 
for $l=1,\ldots,k$.
Assume that $J$ is non periodic.
Then the following holds:
\begin{enumerate}
%(i)
\item
If $J^{'}\in\nset{*}_{A}$ is not in $\{J^{n}:n\geq 1\}$,
then there is $n_{0}\in{\bf N}$ such that
$\pi(s_{J}^{*})^{n}\pi(s_{J^{'}})\Omega=0$ for some $n\geq n_{0}$.
%(ii)
\item
If $v\in {\cal H}$ satisfies $<v|\Omega>=0$, then
$\lim_{n\to\infty}\pi(s_{J}^{*})^{n}v=0$.
\end{enumerate}
\end{lem}
%
% Proof
%
\pr
We simply denote $\pi(s_{i})$ by $s_{i}$ for $i\edot$.

\noindent
(i)
If $J^{'}\in\nset{l}_{A}$ for $1\leq l<k$,
then the non-periodicity of $J$ implies
$s_{J}^{*}s_{J^{'}}\Omega=\delta_{(j_{1},\ldots,j_{l}),J^{'}}\cdot
s_{j_{l+1},\ldots,j_{k}}^{*}\Omega=0$.
If $J^{'}=J_{1}^{'}\cup J_{2}^{'}$ and $|J_{1}^{'}|=nk$
and $|J_{2}^{'}|=l$ for $l=1,\ldots,k-1$,
then $(s_{J}^{*})^{n+1}s_{J^{'}}\Omega=\delta_{J^{n},J^{'}_{1}}
s_{J}^{*}s_{J_{2}^{'}}\Omega=0$ by the previous case.

\noindent
(ii)
By Lemma \ref{lem:cancy}, there is a family
$\{J_{m}^{'}\in\nset{*}_{A}:m\in{\bf N}\}$ such that
$\{s_{J_{m}^{'}}\Omega\}_{m\in {\bf N}}$
is a complete orthonormal basis of ${\cal H}$ and $J_{1}^{'}=J$.
When $<v|\Omega>=0$, we can denote
$v=\sum_{n=2}^{\infty}a_{n}s_{J_{n}^{'}}\Omega$.
If $m\geq 2$, then $J_{m}^{'}\not\in\{J^{n}:n\geq 1\}$.
Therefore there is $n_{0}\in{\bf N}$
such that $(s_{J}^{*})^{n}s_{J_{m}^{'}}\Omega=0$ for $n\geq n_{0}$ by (i).
Hence $\|(s_{J}^{*})^{n}v\|$ is monotone decreasing.
Hence the statement holds.
\qedh

%
% Lemma 4.2
%
\begin{lem}
\label{lem:univec}
Let $({\cal H},\pi)$ be $P(J)$ for $J\in\nset{*}_{A,c}$ 
and $\Omega,\Omega^{'}$ be vectors of ${\cal H}$
such that $\pi(s_{J})\Omega=\Omega$ and $\pi(s_{J})\Omega^{'}=\Omega^{'}$.
If $J$ is non periodic, then $\Omega^{'}=c\Omega$ for some $c\in {\bf C}$.
\end{lem}
%
% Proof
%
\pr
By assumption and Lemma \ref{lem:cancy},
there is a set $\Lambda\subset \nset{*}_{A}$
such that $\Omega^{'}$ is written as 
$c\Omega+\sum_{J^{''}\in\Lambda}a_{J^{''}}\pi(s_{J^{''}})\Omega$
where $<\pi(s_{J^{''}})\Omega|\Omega>=0$ for each $J^{''}\in\Lambda$,
and $\pi(s_{J})^{*}\Omega^{'}=\Omega^{'}$.
By Lemma \ref{lem:lemirre}, $\Omega^{'}=c\Omega$.
\qedh

%
% Lemma 4.3
%
\begin{lem}
\label{lem:irred}
If $J\in\nset{*}_{A,c}$ is non periodic,
then $P(J;z)$ is irreducible for any $z\in U(1)$.
Especially, if $J$ is non periodic, then $P(J)$ is irreducible.
\end{lem}
%
% Proof
%
\pr
Assume that $J$ is non periodic and $({\cal H},\pi)$
is $P(J)$ with the GP vector $\Omega$.
For $v\in {\cal H}$, $v\ne 0$, there is $J^{'}\in\nset{*}_{A}$
such that $<\pi(s_{J^{'}}^{*})v|\Omega>\ne 0$.
Therefore we can always replace
$v$ and $\pi(s_{J^{'}}^{*})v$.
Assume that $v=\Omega+y$ such that $<y|\Omega>=0$.
Then $\lim_{n\to\infty}\pi((s_{J}^{*})^{n})y=0$
by Lemma \ref{lem:lemirre}.
Hence $\lim_{n\to\infty}\pi((s_{J}^{*})^{n})v=\Omega$
and $\Omega\in\overline{\pi(\coa)v}$.
Because $\Omega$ is a cyclic vector, $\overline{\pi(\coa)v}={\cal H}$.
Therefore $({\cal H},\pi)$ is irreducible.
By this and Lemma \ref{lem:guageaction},
$P(J;z)$ is irreducible for each $z\in U(1)$.
\qedh

%
% Lemma 4.4
%
\begin{lem}
\label{lem:vanichain}
Let
$({\cal H},\pi)$
be $P(J)$ with
the GP vector $\Omega$
for
$J=(j_{n})_{n\in{\bf N}}
\in \nset{\infty}_{A}$.
Assume that $J$ is non \evp\ and
define $J_{n}\equiv (j_{1},\ldots,j_{n})$ for $n\in {\bf N}$.
Then the following holds:
\begin{enumerate}
%(i)
\item
For each $m,k\geq 1$ and $J^{'}\in\nset{k}_{A}$
so that $J^{'}\ne J_{k}$, there is $n_{0}\in{\bf N}$
such that $\pi(s_{J_{n}}s_{J_{n}}^{*}s_{J^{'}} 
s_{J_{m}}^{*})\Omega=0$ for each $n\geq n_{0}$.
%(ii)
\item
For each $k\geq 1$ and $J^{'}\in\nset{k}_{A}$
so that $J^{'}\ne J_{k}$, there is $n_{0}\in{\bf N}$
such that $\pi(s_{J_{n}}s_{J_{n}}^{*}s_{J^{'}})\Omega=0$
for each $n\geq n_{0}$.
%(iii)
\item
For each $m\geq 1$, there is $n_{0}\in{\bf N}$
such that $\pi(s_{J_{n}}s_{J_{n}}^{*}s_{J_{m}}^{*})\Omega=0$
for each $n\geq n_{0}$.
%(iv)
\item
If $y\in{\cal H}$ satisfies $<y|\Omega>=0$,
then $\lim_{n\to \infty}\pi(s_{J_{n}}s_{J_{n}}^{*})y=0$.
\end{enumerate}
\end{lem}
%
% Proof
%
\pr
We simply denote $\pi(s_{i})$ by $s_{i}$.

\noindent
(i) and (ii) follow by assumption of $J^{'}$.

\noindent
(iii)
Since $J$ is non \evp, there is $n_{0}\in{\bf N}$
such that $s_{J_{n}}^{*}s_{J_{m}}^{*}s_{J_{n+m}}$
$=0$ for each $n\geq n_{0}$.
Therefore $s_{J_{n}}^{*}s_{J_{m}}^{*}\Omega=s_{J_{n}}^{*}s_{J_{m}}^{*}
s_{J_{n+m}}s_{J_{n+m}}^{*}\Omega=0$ for each $n\geq n_{0}$.
Hence the statement holds.

\noindent
(iv)
Denote $e_{1}\equiv \Omega$ and $e_{m}\equiv s_{J_{m-1}}^{*}\Omega$
for $m\geq 2$.
By Lemma \ref{lem:cancy}, there is a family
$\{K_{n,m}\}_{n,m\in{\bf N}}\subset\nset{*}_{A}$
such that $\{s_{K_{n,m}}e_{m}\}_{n,m\in{\bf N}}$
is a complete orthonormal basis of ${\cal H}$.
If $y\in{\cal H}$ satisfies $<y|\Omega>=0$, then we can denote 
$y=\sum_{m\geq 2}\sum_{n\geq 1}a_{n,m}s_{K_{n,m}}e_{m}$.
Therefore $\|s_{J_{n}}s_{J_{n}}^{*}y\|$ is monotone decreasing
by (i),(ii),(iii). 
Hence the statement holds.
\qedh

%
% Lemma 4.5
%
\begin{lem}
\label{lem:irredd}
If $J\in\nset{\infty}_{A}$ is non \evp, then $P(J)$ is irreducible.
\end{lem}
%
% Proof
%
\pr
Let $({\cal H},\pi)$ be $P(J)$ with the GP vector
$\Omega$ for $J=(j_{n})_{n\in{\bf N}}\in \nset{\infty}_{A}$.
Denote $J_{n}\equiv (j_{1},\ldots,j_{n})$.
Let $v\in {\cal H}$ such that $v\ne 0$.
We can assume that $<v|\Omega>\ne 0$ by replacing
$v$ by $\pi(s_{J^{'}}s_{J^{''}}^{*})v$ for suitable 
$J^{'},J^{''}\in \nset{*}_{A}$ if it is necessary.
We can assume that $v=\Omega+y$ for some $y\in {\cal H}$
such that $<y|\Omega>=0$.
By Lemma \ref{lem:vanichain},
$\Omega=\lim_{n\to\infty}
\pi(s_{J_{n}}s_{J_{n}}^{*})y\in\overline{\pi(\coa)y}$.
Because $\Omega$ is a cyclic vector,
$\overline{\pi(\coa)y}={\cal H}$ and
the statement holds.
\qedh

\noindent
The inverse of Lemma \ref{lem:irredd} is shown in $\S$\ref{section:fifth}.

%%%%%%%%%%%%%%%%%%%%%%%%%%%%%%%%%%%%%%%%%%%
%
% subsection 4.2
%
\ssft{Equivalence}
\label{subsection:fourthtwo}
Recall the equivalence among multiindices
in $\S$\ref{subsection:secondone}.
%
% Lemma 4.6
%
\begin{lem}
\label{lem:equlem}
Let $({\cal H},\pi)$ be $P(J)$ with the GP vector
$\Omega$ for $J\in\nset{*}_{A,c}$.
Assume that $J$ is non periodic
and choose $J^{'}\in\nset{*}$ such that $J^{'}\not\sim J$.
\begin{enumerate}
%(i)
\item
If $\Omega^{'}\in {\cal H}$ satisfies
that $\pi(s_{J^{'}})\Omega^{'}=\Omega^{'}$,
then $<\Omega|\Omega^{'}>=0$.
%(ii)
\item
If $v\in {\cal H}$, then
$\lim_{n\to \infty}\pi(s_{J^{'}}^{*})^{n}v=0$.
\end{enumerate}
\end{lem}
%
% Proof
%
\pr
(i)
Assume that $|J|=k$ and $|J^{'}|=l$.
Because of the non-periodicity of $J$, $J^{l}\ne (J^{'})^{k}$.
This implies the statement.

\noindent (ii)
If $v=\pi(s_{J^{''}})\Omega$, then
$\pi((s_{J^{'}}^{*})^{n_{0}+n})v=\delta_{(J^{'})^{m},J^{''}}
\cdot\pi((s_{J^{'}}^{*})^{n})\Omega\to 0$ 
when
$n\to \infty$
by Lemma \ref{lem:lemirre}
and (i).
Because
any $v\in {\cal H}$
is a limit of linear combination
of
$\{\pi(s_{J^{''}})\Omega:J^{''}\in\nset{*}_{A}\}$,
the statement holds.
\qedh

%
% Lemma 4.7
%
\begin{lem}
\label{lem:equivalence}
Let $J,J^{'}\in\nset{*}_{A,c}$.
\begin{enumerate}
%(i)
\item 
If $J$ and $J^{'}$ are non periodic and $J\not \sim J^{'}$,
then $P(J)\not\sim P(J^{'})$.
%(ii)
\item
For $z,z^{'}\in U(1)$, if $(J,z)\sim (J^{'},z^{'})$,
then $P(J;z)\sim P(J^{'};z^{'})$.
Especially, if $J\sim J^{'}$, then $P(J)\sim P(J^{'})$.
\end{enumerate}
\end{lem}
%
% Proof
%
\pr
(i)
Assume that $P(J)\sim P(J^{'})$.
Then there is a representation $({\cal H},\pi)$ of $\coa$
which is $P(J)$ and $P(J^{'})$.
Assume that $\Omega,\Omega^{'}\in {\cal H}$
are GP vectors with respect to $P(J)$ and $P(J^{'})$, respectively.
Then $<\Omega^{'}|v>=<\Omega^{'}|\pi((s_{J^{'}}^{*})^{n})v>\to 0$
when $n\to \infty$ by Lemma \ref{lem:equlem}.
Hence $\Omega^{'}=0$.
Therefore this is contradiction.
Hence the statement holds.

\noindent
(ii)
Assume that $J\sim J^{'}$ and $J=(j_{1},\ldots,j_{k})\in\nset{k}_{A,c}$.
Let $({\cal H},\pi)$ be $P(J)$ with the GP vector $\Omega$.
Then $\pi(s_{J})\Omega=\Omega$.
By assumption, there is $\sigma\in {\bf Z}_{k}$ such that $J^{'}=\sigma(J)$.
Define $\Omega^{'}\equiv\pi(s_{j_{\sigma(1)}}\cdots s_{j_{\sigma(k)}})\Omega$.
Then $\pi(s_{J^{'}})\Omega^{'}=\Omega^{'}$
and $\Omega^{'}$ is a cyclic vector of $({\cal H},\pi)$.
Hence $P(J)\sim ({\cal H},\pi)\sim P(J^{'})$.
If $(J,z)\sim (J^{'},z^{'})$, then
$J\sim J^{'}$ and $z=z^{'}$ by definition.
Then  $P(J)\sim  P(J^{'})$ by the previous result.
By Lemma \ref{lem:guageaction},
the statement holds.
\qedh

%
% Proposition 4.8
%
\begin{prop}
\label{prop:equivthe}
For $J,J^{'}\in\nset{*}_{A,c}$, assume that $J,J^{'}$ are non periodic.
Then the following holds:
\begin{enumerate}
%(i)
\item
$P(J)\sim P(J^{'})$ if and only if $J\sim J^{'}$.
%(ii)
\item
For $z,z^{'}\in U(1)$, $P(J;z)\sim P(J^{'};z^{'})$
if and only if $(J,z)\sim (J^{'},z^{'})$.
\end{enumerate}
\end{prop}
%
% Proof
%
\pr
(i)
By Lemma \ref{lem:equivalence}, the statement holds.

\noindent
(ii)
If $P(J)\sim P(J^{'};z)$,
then $J\sim J^{'}$ by the proof of Lemma \ref{lem:equivalence}.
If $P(J)$ and  $P(J;z)$ are equivalent, then there is a representation
$({\cal H},\pi)$ which is $P(J)$ and $P(J;z)$
with GP vectors $\Omega$ and $\Omega^{'}$, respectively.
Then $\pi(s_{J})\Omega=\Omega$ and $\pi(s_{J})\Omega^{'}=z\Omega^{'}$.
From this, $z=1$ or $<\Omega|\Omega^{'}>=0$.
If $<\Omega|\Omega^{'}>=0$, then
$\pi((s_{J}^{*})^{n})\Omega^{'}\to 0$ when $n\to \infty$.
This is contradiction.
Therefore $P(J)\not\sim P(J;z)$ when $z\ne 1$.
On the other hand, $P(J;z)\sim P(J^{'};z^{'})$
if and only if $P(J)\sim P(J^{'};z^{'})\circ
\gamma_{\bar{z}^{1/k}}=P(J^{'};z^{'}\bar{z}^{l/k})$
when $|J|=k$ and $|J^{'}|=l$
by Lemma \ref{lem:guageaction}.
Therefore $P(J;z)\sim P(J^{'};z^{'})$
if and only if $J\sim J^{'}$ and $z^{'}\bar{z}=1$.
This implies the statement.
\qedh

\noindent
In $\S$\ref{section:fifth}, we show the statement in
Proposition \ref{prop:equivthe} without assumption of non periodicity.

%
% Lemma 4.9
%
\begin{lem}
\label{lem:chaineq}
For $J,J^{'}\in\nset{\infty}_{A}$, 
if $J\not \sim J^{'}$, then $P(J)\not\sim P(J^{'})$.
\end{lem}
%
% Proof.
%
\pr
Assume that $J=(j_{n})_{n\in {\bf N}}$,
$J^{'}=(j_{n}^{'})_{n\in {\bf N}}$, $J\not\sim J^{'}$
and $P(J)\sim P(J^{'})$.
Then there is a representation
$({\cal H},\pi)$ which is
$P(J)$ and $P(J^{'})$ with GP vectors
$\Omega$ and $\Omega^{'}$, respectively.
Define $e_{n}\equiv \pi(s_{j_{n}}^{*}\cdots
s_{j_{1}}^{*})\Omega$
and $e^{'}_{n}\equiv \pi(s_{j_{n}^{'}}^{*}\cdots
s_{j_{1}^{'}}^{*})\Omega^{'}$
where $J_{n}\equiv (j_{1},\ldots,j_{n})$
and $J_{n}^{'}\equiv (j_{1}^{'},\ldots,j_{n}^{'})$
for $n\geq 1$.
Then
$<\Omega|\Omega^{'}>=<\pi(s_{J_{n}})
e_{n+1}|\pi(s_{J_{n}^{'}})e_{n+1}^{'}>=0$
for some $n\in{\bf N}$
because $J\not\sim J^{'}$.
In the same way, we see that $<e_{m}|e_{n}^{'}>=0$
for each $n,m\in {\bf N}$.
From this, $<\pi(s_{K})e_{m}|e_{n}^{'}>=0$
for each $K\in \nset{*}_{A}$ and $n,m\in {\bf N}$.
Therefore $<v|e_{n}^{'}>=0$ for each $v\in {\cal H}$
and $n\in {\bf N}$.
Hence $e_{n}^{'}=0$ for each $n\in{\bf N}$.
This contradicts with the choice of
$\{e_{n}^{'}\}_{n\in{\bf N}}$.
Therefore $P(J)\not \sim P(J^{'})$.
\qedh

%
% Proposition 4.10
%
\begin{prop}
\label{prop:equcha}
For $J,J^{'}\in\nset{\infty}_{A}$,
$P(J)\sim P(J^{'})$ if and only if $J\sim J^{'}$.
\end{prop}
%
% Proof
%
\pr
By Lemma \ref{lem:chaineq}, it is sufficient to
show that $J\sim J^{'}$ implies that $P(J)\sim P(J^{'})$.
Assume that $J=(j_{n})_{n\in{\bf N}}$,
$J^{'}=(j_{n}^{'})_{n\in{\bf N}}$ and $J\sim J^{'}$.
Then there are $p\in {\bf Z}$ and $M\geq 1$
such that $j_{n}^{'}=j_{n+p}$ for each $n\geq M$.
If $({\cal H},\pi)$ is $P(J)$ with the GP vector $\Omega$, then
define $t_{n}\equiv s_{j_{1}}\cdots s_{j_{n}}$
and  $t_{n}^{'}\equiv s_{j_{1}^{'}}\cdots s_{j_{n}^{'}}$ for 
$n\in {\bf N}$ and $\Omega^{'}\equiv \pi(t_{M}^{'}t_{M+p}^{*})\Omega$.
Then $\pi((t^{'}_{M+n})^{*})\Omega^{'}=\pi(t_{M+n+p}^{*})\Omega$
for each $n\geq 1$.
Therefore $\{\pi((t_{l}^{'})^{*})\Omega^{'}\}_{l\geq 1}$
is a chain of $\pi$ by $J^{'}$.
Because $\Omega$ is a cyclic vector, $\Omega^{'}$ is.
Hence the statement holds.
\qedh

%%%%%%%%%%%%%%%%%%%%%%%%%%%%%%%%%%%%%%%%%%%%%%%
%
% Section 5
%
\sftt{Decomposition}
\label{section:fifth}
%
% Proposition 5.1
%
\begin{prop}
\label{prop:decc}
For $(J,c)\in\nset{*}_{A,c}\times U(1)$ and $p\geq 1$,
%
% Equation 5.1
%
\begin{equation}
\label{eqn:deg}
P(J^{p};c)\sim\bigoplus_{j=1}^{p}P(J;c^{1/p}\xi^{j})
\end{equation}
where $\xi\equiv e^{2\pi\sqrt{-1}/p}$.
Especially,
%
% Equation 5.2 
%
\begin{equation}
\label{eqn:degp}
P(J^{p})\sim\bigoplus_{j=1}^{p}P(J;\xi^{j}).
\end{equation}
\end{prop}
%
% Proof
%
\pr
Assume that $J=(j_{l})_{l=1}^{k}\in\nset{k}_{A,c}$.
Let $({\cal H}_{j},\pi_{j})$ be $P(J;\xi^{j})$ with 
the GP vector $\Omega_{j}$ for $j=1,\ldots,p$.
Let $\Omega\equiv p^{-1/2}\sum_{j=1}^{p}\Omega_{j}
\in{\cal H}\equiv {\cal H}_{1}\oplus\cdots\oplus {\cal H}_{p}$
and $\pi\equiv \pi_{1}\oplus\cdots\oplus\pi_{p}$.
Then $\pi(s_{J^{p}})\Omega=\Omega$ and
$\{\pi(s_{j_{l}}\cdots s_{j_{k}}s_{J^{a}})
\Omega:a=0,\ldots,p-1,\,
l=1,\ldots,k\}$ is an orthonormal family.
Therefore $V\equiv \overline{\pi(\coa)\Omega}$
is a $P(J^{p})$-component of ${\cal H}$.
On the other hand,
${\rm Lin}\langle\{\pi(s_{J^{q}})\Omega:q=1,\ldots,p\}\rangle
={\rm Lin}\langle\{\Omega_{q}:q=1,\ldots,p\}\rangle$.
Hence $\Omega_{i}\in V$ and
${\cal H}_{i}\subset V$ for $i=1,\ldots,p$.
Therefore ${\cal H}=V$
and  $({\cal H},\pi)$ is $P(J^{p})$.
From this, we obtain (\ref{eqn:degp}).
By Lemma \ref{lem:guageaction} and (\ref{eqn:degp}),
(\ref{eqn:deg}) is verified.
\qedh
%
% Corollary 5.2
%
\begin{cor}
\label{cor:decc}
\begin{enumerate}
%(i)
\item
For $(J,z)\in\nset{*}_{A,c}\times U(1)$, 
$P(J;z)$ is irreducible if and only if $J$ is non periodic.
Especially, for $J\in\nset{*}_{A,c}$,
$P(J)$ is irreducible if and only if $J$ is non periodic. 
%(ii)
\item
For $J,J^{'}\in\nset{*}_{A,c}$ and $z,z^{'}\in U(1)$,
$P(J;z)\sim P(J^{'};z^{'})$ if and only if $(J,z)\sim (J^{'},z^{'})$.
%(iii)
\item
The decomposition in (\ref{eqn:deg}) is multiplicity free
when $J$ is non periodic.
\end{enumerate}
\end{cor}
%
% Proof
%
\pr
(i)
By Lemma \ref{lem:irred} and Theorem \ref{cor:decc}, the statement holds.

\noindent
(ii)
If $(J,z)\sim (J^{'},z^{'})$, then $P(J;z)\sim P(J^{'};z^{'})$
by Lemma \ref{lem:equivalence}.
Assume that $P(J;z)\sim P(J^{'};z^{'})$.
If $J$ and $J^{'}$ are non periodic,
then the statement is shown in
Proposition \ref{prop:equivthe}.
If $J$ is periodic, then $P(J;z)$ is not irreducible
by (i) and decomposed into direct sum of finite irreducible components
by Theorem \ref{cor:decc}.
Therefore $P(J^{'};z^{'})$ must not be irreducible.
By (i), $J^{'}$ is periodic.
By comparing irreducible components
of $P(J;z)$ and those of $P(J^{'};z^{'})$ and Proposition \ref{prop:equivthe},
we see that their sets of irreducible components
coincide up to unitary equivalence.
From this, $(J,z)\sim (J^{'},z^{'})$.

\noindent
(iii) By (i) and (ii), the assertion holds.
\qedh

\noindent
By Proposition \ref{prop:equcha}
and Corollary \ref{cor:decc}, we have the following:
%
% Proposition 5.3
%
\begin{prop}
\label{prop:cycheq}
For $J,J^{'}\in\nset{\#}_{A,c}$, $P(J)\sim P(J^{'})$ if and only if $J\sim J^{'}$.
\end{prop}

%%%%%%%%%%%%%%%%%%%%%%%%%%%%%%%%%%%%%%%%%%%%%%%%%%%%%%%%%%%%
%
% Lemma 5.4
% 
\begin{lem}
\label{lem:pchain}
For $J\in \nset{*}_{A,c}$, $P(J^{\infty})\sim \int_{U(1)}^{\oplus}P(J;c)\,d\eta(c)$
where $\eta$ is the Haar measure of $U(1)$.
\end{lem}
%
% Proof
%
\pr
Recall ${\cal T}(A;j)$ in (\ref{eqn:tree}).
For $J=(j_{n})_{n=1}^{k}\in \nset{k}_{A,c}$, define
$\Lambda^{'}\equiv\Lambda^{'}_{1}\sqcup\Lambda^{'}_{2}\sqcup\Lambda^{'}_{3}$
where $\Lambda_{1}^{'}\equiv{\bf Z}$,
$\Lambda_{2}^{'}\equiv\coprod_{n\in{\bf Z}}\Lambda_{2,n}^{'}$,
$\Lambda_{2,n}^{'}\equiv\{(j,n):j\in {\cal T}^{(1)}(A;j_{n}),\,j\ne j_{n-1}\}$ 
and $j_{kn+m}\equiv j_{m}$ for each $n\in {\bf Z}$ and $m=1,\ldots,k$
and $\Lambda_{3}^{'}\equiv\coprod_{(j,n)\in
\Lambda_{2}^{'}}{\cal T}(A;j)\times \{(j,n)\}$.
Define a subset of $\Lambda^{'}$ by
\[D(g_{i})\equiv\{kn+m\in{\bf Z}:n\in {\bf Z},\, a_{ij_{m}}=1\}\sqcup
\{({\cal T}(i;A)\times {\bf Z})\cap \Lambda^{'}\}\]
and let $g=\{g_{i}\}_{i=1}^{N}$ be a family of maps by
$g_{i}:D(g_{i})\to \Lambda^{'}$;
\[
\left\{
\begin{array}{l}
g_{i}(kn+m)\equiv 
\left\{
\begin{array}{ll}
kn+m-1\quad&(i=j_{m-1}),\\
&\\
(i,kn+m)\quad&(\mbox{otherwise})\\
\end{array}
\right.\quad
(kn+m\in\Lambda_{1}^{'}(J)\cap D(g_{i}))\\
\\
g_{i}(J^{'},n)\equiv(\{i\}\cup J^{'},\, n)\quad\quad\quad
((J^{'},n)\in(\Lambda_{2}^{'}\sqcup\Lambda_{3}^{'})\cap D(g_{i})).
\end{array}
\right.
\]
Then $g$ is an $A$-\bfs\ on $\Lambda^{'}$ and $g$ is $P(J^{\infty})$
and $g_{J}(kn+1)=k(n-1)+1$.
Hence $(l_{2}(\Lambda^{'}),\pi_{g})$ is $P(J^{\infty})$.

Let $f$ and $\Lambda(A,J)$ be in Lemma \ref{lem:existcycle}.
Define 
${\cal K}\equiv l_{2}(\Lambda(A,J))$ and	
a unitary $T$ from
$l_{2}(\Lambda^{'})$ to $L_{2}(U(1),{\cal K})$ by
\[Te_{kn+m}\equiv \zeta_{n}\otimes e_{J_{m}},\quad
Te_{j,kn+m}\equiv \zeta_{n}\otimes e_{j,J_{m}},\quad
Te_{x,j,kn+m}\equiv \zeta_{n}\otimes e_{x,j,J_{m}}\]
for $n\in {\bf Z}$ and $m=1,\ldots,k$.
Let $\Pi\equiv {\rm Ad}T\circ \pi_{g}$.
Then
$\Pi(s_{J})(\zeta_{n}\otimes e_{J_{1}})=\zeta_{n-1}\otimes e_{J_{1}}$.
For $c\in U(1)$, define a representation $\pi_{c}$ of $\coa$ on ${\cal K}$	
by
\[\pi_{c}(s_{i})e_{x,J_{l}}\equiv \chi_{D(f_{i})}(x,J_{l-1})
\cdot e_{f_{i}(x,J_{l-1})},\quad \!\!
\pi_{c}(s_{i})e_{x,J_{1}}\equiv \bar{c}\cdot \chi_{D(f_{i})}(x,J_{1})
\cdot e_{f_{i}(x,J_{1})}\]
for $l=2,\ldots,k$, $i\edot$
and $x\in \nset{*}_{A}$ such that $(x,J_{l})\in D(f_{i})$
where $(0,J_{l})\equiv J_{l}$.
By Lemma \ref{lem:existcycle},
we can verify that $({\cal K},\pi_{c})$ is $P(J;\bar{c})$.
Further $(\Pi(s_{i})\phi)(c)=\pi_{c}(s_{i})\phi(c)$
for each $i\edot$, $\phi\in L_{2}(U(1),{\cal K})$ and $c\in U(1)$.
Therefore
$P(J^{\infty})\sim \Pi\sim \int_{U(1)}^{\oplus}\pi_{c}\,d\eta(c)
\sim \int_{U(1)}^{\oplus}P(J;\bar{c})\,d\eta(c)
\sim \int_{U(1)}^{\oplus}P(J;c)\,d\eta(c)$.
\qedh

%%%%%%%%%%%%%%%%%%%%%%%%%%%%%%%%%%%%%%%%%%%%%%%%%%%%%%%%%%%%

%
% Proposition 5.5
%
\begin{prop}
\label{prop:decochain}
If $K\in \nset{\infty}_{A}$ is \evp, 
then there uniquely exists $J\in [1,\ldots,N]_{A}^{*}$ such that
%
% Equation 5.3
%
\begin{equation}
\label{eqn:deceq}
P(K)\sim \int_{U(1)}^{\oplus}P(J;c)\, d\eta(c).
\end{equation}
\end{prop}
%
% Proof
%
\pr
By the assumption of $K$, there is $J\in[1,\ldots,N]_{A}^{*}$
such that $P(K)\sim P(J^{\infty})$.
By Lemma \ref{lem:pchain}, the existence is verified.
By assumption, $P((J^{'})^{\infty})\sim P(K)\sim P(J^{\infty})$.
This is equivalent that $(J^{'})^{\infty}\sim J^{\infty}$.
Since $J,J^{'}\in[1,\ldots,N]_{A}^{*}$,
the uniqueness holds.
\qedh

%
% Proposition 5.6
%
\begin{prop}
\label{prop:irretwo}
For $K\in \nset{\infty}_{A}$, $P(K)$ is irreducible if and only if
$K$ is non \evp.
\end{prop}
%
% Proof.
%
\pr
By Lemma \ref{lem:irredd} and Proposition \ref{prop:decochain},
the statement holds.
\qedh

%%%%%%%%%%%%%%%%%%%%%%%%%%%%%%%%%%%%%%%%%%%
%
% Proposition 5.7
%
\begin{prop}
\label{prop:deco}
For $A\in\mnz$,
let $({\cal H},\pi)$ be a permutative
representation of $\coa$,
and $<1,\ldots,N>^{\#}_{A}$ and $[1,\ldots,N]^{\#}_{A}$
be in (\ref{eqn:alldef}).
\begin{enumerate}
%(i)
\item
The following decomposition into cyclic subspaces holds:
%
% Equation 5.4
%
\begin{equation}
\label{eqn:deco}
({\cal H},\pi)\quad
\sim \bigoplus_{J\in
<1,\ldots,N>^{\#}_{A}
}P(J)^{\oplus \nu_{J}}
\end{equation}
where $\nu_{J}$ is the multiplicity of $P(J)$.
Furthermore (\ref{eqn:deco}) is unique up to unitary equivalence.
%(ii)
\item
The following irreducible decomposition holds:
%%%%%%%%%%%%%%%%%%%%%%%%%%%%%%%%%%%%
\[{\cal H}\quad=\bigoplus_{J\in[1,\ldots,N]_{A}^{*}}{\cal H}_{J}
\quad \oplus\bigoplus_{K\in[1,\ldots,N]_{A}^{\infty}}{\cal H}_{K},\]
\[
\!\!\!\!{\cal H}_{J}=\bigoplus_{p\geq 1}
\left\{\bigoplus_{j=1}^{p}{\cal H}_{J,p,j}\right \}
^{\oplus\nu_{J,p}}\!\!\oplus\left\{\int^{\oplus}_{U(1)}
{\cal H}_{J,\infty,z}\,dm(z)
\right\}^{\oplus \nu_{J,\infty}}\!\!\!,\quad
{\cal H}_{K}=\{{\cal H}_{K,0}\}^{\oplus \nu_{K}}\]
where ${\cal H}_{J,p,j}$,
${\cal H}_{J,\infty,z}$ and ${\cal H}_{K,0}$ are
$P(J;e^{2\pi\sqrt{-1}j/p})$, $P(J;z)$ and  $P(K)$, respectively,
and $\nu_{J,p}$ and $\nu_{K}$ are multiplicities.
\end{enumerate}
\end{prop}
%
% Proof
%
\pr
(i)
By Lemma \ref{lem:decothm} and
Proposition \ref{prop:cycheq}, the statement holds.

\noindent
(ii)
Proposition \ref{prop:decc} and Proposition \ref{prop:decochain}
imply the decomposition.
\qedh

Assume that there are two irreducible decompositions
of a given permutative representation $({\cal H},\pi)$ of $\coa$.
If there is no direct integral component, then the uniqueness follows.
If there is a direct integral decomposition on $U(1)$ 
as a style in (\ref{eqn:deceq}),
then the uniqueness holds in a sense of Proposition \ref{prop:decochain}.
In consequence, the irreducible decomposition of permutative representation
as a form in Proposition \ref{prop:deco} (ii) is unique up to unitary equivalence.

%
% Proposition 5.8
%
\begin{prop}
\label{prop:comp}
For any $A\in\mnz$, any permutative representation
of $\coa$ is completely reducible and irreducible decomposition 
as a form in Proposition \ref{prop:deco} (ii)
is unique up to unitary equivalence.
\end{prop}

%
% Proposition 5.9
%
\begin{prop}
\label{prop:matwoo}
For $A\in\mnz$, let $({\cal H},\pi)$
be a permutative representation of $\coa$ with phases.
Then the following unique decomposition into cyclic representations 
up to unitary equivalence holds:
\[({\cal H},\pi)\sim \bigoplus_{(J,c)\in[1,\ldots,N]_{A}^{*}\times U(1)}
P(J;c)^{\oplus \nu_{J,c}}\oplus\bigoplus_{K\in<1,\ldots,N>_{A}^{\infty}}
P(K)^{\oplus \nu_{K}}
\]
where $\nu_{J,c}$ and $\nu_{K}$ are multiplicities.
Especially, if  $({\cal H},\pi)$ is cyclic, then $({\cal H},\pi)$
is equivalent to either $P(J;c)$ or $P(K)$.
\end{prop}
%
% Proof
%
\pr
By assumption, there is a complete orthonormal basis
$\{e_{n}\}_{n\in\Lambda}$ of ${\cal H}$,
$\{\Lambda_{i}\}_{i=1}^{N}$
and $\{(z_{i,n},m_{i,n})\in U(1)
\times \Lambda:(i,n)
\in\nset{}\times \Lambda_{i}\}$ such that
$\pi(s_{i})e_{n}=z_{i,n}\chi_{\Lambda_{i}}(n)\cdot e_{m_{i,n}}$
for  each $(i,n)\in\nset{}\times \Lambda$.
Define a new permutative representation $({\cal H},\pi_{0})$
of $\coa$ by
$\pi_{0}(s_{i})e_{n}\equiv \chi_{\Lambda_{i}}(n)\cdot e_{m_{i,n}}$
for $(i,n)\in\nset{}\times \Lambda$.
By Proposition \ref{prop:deco} (i), 
$\pi_{0}$ is decomposed into the direct sum of permutative representations:
\[
\pi_{0}\sim \bigoplus_{
J\in<1,\ldots,N>^{*}_{A}}
P(J)^{\oplus \nu_{J}}\oplus
\bigoplus_{
K\in<1,\ldots,N>^{\infty}_{A}}
P(K)^{\oplus \nu_{K}}.
\]
Therefore $\pi_{0}|_{V}$ is $P(J)$ or $P(K)$
for some subspace $V\subset {\cal H}$.
If $\pi_{0}|_{V}$ is $P(J)$, then there is a cyclic unit vector
$\Omega \in V$ such that $\pi_{0}(s_{J})\Omega=\Omega$.
By definition of $\pi_{0}$, there is $c_{J}\in U(1)$
such that $\pi(s_{J})\Omega=c_{J}\Omega$.
Because $(V,\pi_{0}|_{V})$ is cyclic, $(V,\pi|_{V})$ is.
Therefore $\pi|_{V}$ is $P(J;c_{J})$.
If $\pi_{0}|_{V}$ is $P(K)$, then we see that
$\pi|_{V}$ is $P(K)$ by checking the condition of chain.
In consequence
\[\pi\sim \bigoplus_{J\in<1,\ldots,N>^{*}_{A}}
P(J;c_{J})^{\oplus \nu_{J}}\oplus
\bigoplus_{K\in<1,\ldots,N>^{\infty}_{A}}P(K)^{\oplus \nu_{K}}.\]
When $J$ is periodic, $P(J;c_{J})$ is decomposed 
into the direct sum of elements in
$\{P(J^{'};c^{'}):(J^{'},c^{'})\in [1,\ldots,N]^{*}_{A}\times U(1)\}$
by Theorem \ref{cor:decc}. 
Hence the statement holds.
\qedh

\noindent
{\it Proof of Theorem \ref{Thm:matwo}.}
(i) and (ii) follow from  Proposition \ref{prop:matwoo}.
(iii) follows from Corollary \ref{cor:decc} (i) and Proposition \ref{prop:irretwo}.
(iv), (v) and (vi) follow from 
Corollary \ref{cor:decc} (ii) and Proposition \ref{prop:cycheq}.
\qedh

%%%%%%%%%%%%%%%%%%%%%%%%%%%%%%%%%%%%%%
%
% Section 6
%
\sftt{States and spectrums}
\label{section:sixth}
Fix $A\in\mnz$.
Operator algebraists prefer {\it states} than representations.
Therefore we show states of the \cka s associated with
permutative representations.
%
% Proposition 6.1
%
\begin{prop}
\label{prop:stacycle}
Let $({\cal H},\pi)$ be $P(J)$ with the GP vector $\Omega$
for $J\in\nset{\#}_{A,c}$.
Define a state $\omega$ of $\coa$ by 
$\omega\equiv <\Omega|\pi(\cdot)\Omega>$.
Then the following holds:
\begin{enumerate}
%(i)
\item
When $J=(j_{1},\ldots,j_{k})\in\nset{k}_{A,c}$,
%
% Equation 6.1
%
\begin{equation}
\label{eqn:state}
\omega(s_{J^{'}}s_{J^{''}}^{*})=
\left\{
\begin{array}{ll}
1 \quad &(0\leq\exists p\leq k-1,\,
s.t. \, \,J^{'},J^{''}\in
{\cal I}_{p}(J)),\\
&\\
0\quad&(\mbox{otherwise})\\
\end{array}
\right.
\end{equation}
where ${\cal I}_{p}(J)\equiv \{J^{a}\cup (j_{1},\ldots,j_{p})
\in\nset{*}_{A}:a\geq 0\}$.
%(ii)
\item
When $J=(j_{n})_{n\in{\bf N}}\in\nset{\infty}_{A}$,
%
% Equation 6.2
%
\begin{equation}
\label{eqn:statechain}
\omega(s_{J^{'}}s_{J^{''}}^{*})
=
\left\{
\begin{array}{ll}
1\quad &(\exists m\in {\bf N}\, s.t.\,
J^{'}=J^{''}=(j_{1},\ldots,j_{m})),\\
&\\
0\quad&(\mbox{otherwise}).
\end{array}
\right.
\end{equation}
%(iii)
\item
The GNS representation
of $\coa$
by a state $\omega$
which satisfies
(\ref{eqn:state})
or
(\ref{eqn:statechain})
is equivalent to
$P(J)$.
%(iv)
\item
$\omega$
is pure 
 if and only if
$J$ is non periodic
or
non \evp.
\end{enumerate}
\end{prop}
%
% Proof
%
\pr
(i) and (ii) are verified by direct computation.

\noindent
(iii)
The statement follows from the uniqueness of the GNS representation.

\noindent
(iv)
Corollary \ref{cor:decc} and Proposition \ref{prop:irretwo}
imply the assertion.
\qedh

We consider the spectrum of $\coa$ associated with
permutative representations of $\coa$.
${\rm Spec}\coa$ is the {\it spectrum} of $\coa$
which consists of all unitary equivalence classes
of irreducible representations of $\coa$.
Define ${\rm PSpec}\coa$ by the subset of ${\rm Spec}\coa$
which consists of all unitary equivalence classes
of irreducible permutative representations.
By Proposition \ref{prop:uniqeness}, Corollary \ref{cor:decc}
and Proposition \ref{prop:irretwo},
the following one-to-one correspondence holds:
\[{\rm PSpec}\coa\cong [1,\ldots,N]^{\#}_{A}.\]
By regarding phase factor,
$\{[1,\ldots,N]^{*}_{A}\times U(1)\}\sqcup
[1,\ldots,N]^{\infty}_{A}$
is identified with a subset of ${\rm Spec}\coa$.
Let $A_{1}\equiv \mattwo{1}{1}{0}{1}$ and $A_{2}\equiv \mattwo{1}{1}{1}{0}$.
We see that $\{1,2\}^{*}_{A_{1},c}=\{(1)^{n},(2)^{n}:n\geq 1\}$
and $\{1,2\}^{\infty}_{A_{1}}=\{(1)^{\infty},
(1)^{n-1}\cup(2)^{\infty}:n\geq 1\}$.
Hence $[1,2]^{*}_{A_{1}}=\{(1),(2)\}$, $[1,2]^{\infty}_{A_{1}}=\emptyset$
and $\#{\rm PSpec}\co{A_{1}}=2$.
$\{(1)^{n}\cup (2):n\geq 1\}$
is a proper subset of $[1,2]^{*}_{A_{2}}$.
From this, $\#{\rm PSpec}\co{A_{2}}=\infty$.

%%%%%%%%%%%%%%%%%%%%%%%%%%%%%%%%%%%%%%%%
%
% Section 7
%
\sftt{Examples}
\label{section:seventh}
%%%%%%%%%%%%%%%%%%%%%%%%%%%%%%%%%%%%%%
%
% subsection 7.1
%
\ssft{Naive examples}
\label{subsection:seventhone}
Let $A_{1}\equiv\matthr
{0}{0}{1}
{1}{0}{1}
{1}{1}{1}$ and
$A_{2}\equiv
\matthr
{0}{1}{1}
{1}{0}{1}
{1}{1}{1}$.
Define a representation $(\ltn,\pi)$ of $\co{A_{1}}$ by
\[
\begin{array}{c}
\pi(s_{1})e_{4(n-1)+i}\equiv \delta_{2,i}e_{4(n-1)+1},\quad
\pi(s_{3})e_{n}\equiv e_{4(n-1)+2},\\
\\
\pi(s_{2})e_{4(n-1)+i}\equiv \delta_{1,i}e_{4(n-1)+4}+
\delta_{2,i}e_{4(n-1)+3}\\
\end{array}
\]
for $n\in {\bf N}$ and $i=1,2,3,4$.
Then $(\ltn,\pi)$ is $P(13)$.\\
%
% Proof
%
\pr
Define
$D_{1}\equiv \{4(n-1)+2:n\in {\bf N}\}$,
$D_{2}\equiv \{4(n-1)+1,4(n-1)+2:n\in {\bf N}\}$,
$D_{3}\equiv {\bf N}$ and $f=\{f_{1},f_{2},f_{3}\}$
by $\pi(s_{i})e_{n}\equiv e_{f_{i}(n)}$ for $i=1,2,3$
and $n\in D_{i}$.
Then $f\in\bfset{A_{1}}{{\bf N}}$ and $f$ is cyclic.
We see that $f$ has neither chain nor cycle in $\{n\in{\bf N}:n\geq 3\}$
and $(f_{1}\circ f_{3})(1)=1$.
This implies that $\pi(s_{1}s_{3})e_{1}=e_{1}$.
Hence the statement holds.
\qedh

Define a representation $(l_{2}({\bf N}\times \{1,2\}),\pi)$
of $\co{A_{2}}$ by
\[
\pi(s_{1})e_{n,i}\equiv 
\delta_{2,i}e_{n,1},\quad
\pi(s_{3})e_{n,i} \equiv \delta_{1,i}e_{5(n-1)+4,2}
+\delta_{2,i}e_{5(n-1)+5,2},
\]
\[
\pi(s_{2})e_{5(n-1)+m,i} \equiv 
\delta_{1,i}
e_{5(5(n-1)+m-1)+1,2}
+\delta_{2,i}(\delta_{4,m}e_{5(n-1)+2,2}
+\delta_{5,m}e_{5(n-1)+3,2})
\]
for $i=1,2$, $m=1,\ldots,5$ and $n\in {\bf N}$
where $\{e_{n,i}\}$ is the canonical basis of $l_{2}({\bf N}\times \{1,2\})$.
Then $(l_{2}({\bf N}\times \{1,2\}),\pi)$ is $P(12)$.\\
%
% Proof
%
\pr
We see that $\pi$ is cyclic, $\pi(s_{1}s_{2})e_{1,1}=e_{1,1}$
and there is no cycle except this.
Therefore the statement holds.
\qedh

%%%%%%%%%%%%%%%%%%%%%%%%%%%%%%%%%%%%%%%%%%%%%%%%
%
% Subsection 7.2
%
\ssft{Standard representation}
\label{subsection:seventhtwo}
We define the standard $A$-\bfs\ and the standard 
representation of $\coa$ for a given $A\in\mnz$.

For $A=(a_{ij})\in \mnz$, a data $\mqb$ is called the (canonical) 
{\it $A$-coordinate} if $B_{i}\equiv \left\{\,j\in\{1,\ldots,N\}:a_{ij}=1\,\right\}$,
$M_{i}\equiv a_{i1}+\cdots+a_{iN}$,
$q_{i}:B_{i}\to\{1,\ldots,M_{i}\};
q_{i}(j)\equiv\#\{k\in B_{i}:k\leq j\}$ for $i\edot$.
An $A$-\bfs\ $f^{(A)}=\{f_{i}^{(A)}\}_{i=1}^{N}$
on ${\bf N}$ defined by
\[f_{i}^{(A)}(N(m-1)+j)\equiv N(M_{i}(m-1)+q_{i}(j)-1)+i\quad
(m\in {\bf N},\,j\in B_{i}),\]
\[R(f_{i}^{(A)})\equiv \{N(n-1)+i:n\in{\bf N}\},\,\,
\mbox{$D(f_{i}^{(A)})\equiv \coprod_{j\in B_{i}}
R(f_{j}^{(A)})\,\,(i\edot)$}\]
is called the {\it standard $A$-\bfs}.
We call $(\ltn,\pi_{f^{(A)}})$ by
the {\it standard representation} of $\coa$
and denote $(\ltn,\pi_{S}^{(A)})$.

In order to show the decomposition formula of the standard representation,
we define cycles arising from some finite dynamical system associated with $A$.
For $A=(a_{ij})\in\mnz$, define a map $\varphi_{A}$ on $\nset{}$ by 
% 
% Equation 7.1
%
\begin{equation}
\label{eqn:varphia}
\varphi_{A}(i)\equiv \min\{j\in\nset{}:a_{ij}=1\}\quad(i\edot).
\end{equation}
Then $\nset{}$ contains cycles by $\varphi_{A}$,
that is, $C=\{n_{i}\in \nset{}:i=1,\ldots,m\}$
is a {\it cycle} in $\nset{}$ by $\varphi_{A}$
if $\varphi_{A}(n_{i})=n_{i+1}$ for $i=1,\ldots,m-1$
and $\varphi_{A}(n_{m})=n_{1}$.
%
% Definition 7.1
%
\begin{defi}
\label{defi:cycadef}
\begin{enumerate}
%(i)
\item
$\{C_{i}\}_{i=1}^{k}$ is the $A$-cycle set
if $\{C_{i}\}_{i=1}^{k}$ is the set of all cycles 
in $\nset{}$ by $\varphi_{A}$.
Define $m_{i}\equiv \#C_{i}$ and
$j_{i,1}\equiv \min C_{i}$ for $i=1,\ldots,k$.
%(ii)
\item
${\cal J}_{A}\equiv \{J_{i}\}_{i=1}^{k}\subset
\nset{*}$ is  the $A$-cyclic index set
if $J_{i}\equiv (j_{i,c})_{c=1}^{m_{i}}\in\nset{m_{i}}$ 
and $j_{i,c}\equiv \varphi_{A}^{c-1}(j_{i,1})$
for $c=1,\ldots,m_{i}$
where $j_{i,1}\in\nset{}$ and $m_{i}$ are in (i).
Let
${\cal J}_{A,\infty}\equiv \{(j_{i})_{i=1}^{m}
\in {\cal J}_{A}:\forall l, 
a_{j_{1},l}=\delta_{l,j_{2}},\ldots,
a_{j_{m-1},l}=\delta_{l,j_{m}},a_{j_{m},l}=\delta_{l,j_{1}}\}$,
${\cal J}_{A,1}\equiv {\cal J}_{A}\setminus
{\cal J}_{A,\infty}$.
\end{enumerate}
\end{defi}

%
% Lemma 7.2
% 
\begin{lem}
\label{lem:abfsth}
\begin{enumerate}
%(i)
\item
%\label{lem:strstan}
Let $A\in\mnz$ with the $A$-coordinate $\mqb$.
For $J=(j_{1},\ldots,j_{k})\in\nset{k}_{A,c}$,
the standard $A$-\bfs\ has a $P(J)$-component
if and only if $q_{j_{i}}(j_{i+1})=1$
for each $i=1,\ldots,k-1$ and $q_{j_{k}}(j_{1})=1$.
%(ii)
\item
%\label{lem:nochain}
For any $A\in\mnz$, the standard $A$-\bfs\ has no chain.
%(iii)
\item
For $A\in\mnz$, if ${\cal J}_{A}$ is the $A$-cyclic index set
and ${\cal J}_{A,1}$, ${\cal J}_{A,\infty}$ are in Definition
\ref{defi:cycadef},
then the standard $A$-\bfs\ is decomposed as
$\bigoplus_{J\in{\cal J}_{A,1}}P(J)\oplus
\left\{\bigoplus_{K\in{\cal J}_{A,\infty}}P(K)\right\}^{\oplus\infty}$.
\end{enumerate}
\end{lem}
%
% Proof
%
\pr
(i) and (ii) are verified by their definition directly.

\noindent
(iii)
We denote $f^{(A)}$ by $f$ simply.
Let $f$ be the standard $A$-\bfs.
By Lemma \ref{lem:decothm} and (ii),
$f$ is decomposed into only cycles.
On the other hand, any cycle component of $f$ is one of 
$\{P(J):J\in{\cal J}_{A}\}$ by (i).
Therefore $f$ is decomposed as a direct sum of 
$\{P(J):J\in{\cal J}_{A}\}$ with multiplicities.	
If $J\in{\cal J}_{A,\infty}$,
then $f$ has a $P(J)^{\oplus \infty}$-component.
If $J\in{\cal J}_{A,1}$, then the cycle of $f$ by $J$
is a subset of $\nset{}$.
By definition of ${\cal J}_{A}$,
$P(J)$ appears in $\nset{}$ once for all.
In consequence, the statement holds.
\qedh

%
% Proposition 7.3
%
\begin{prop}
\label{prop:stast}
For $A\in\mnz$, let $(\ltn,\pi_{S}^{(A)})$ 
be the standard representation of $\coa$.
Then the following holds:
\begin{enumerate}
%(i)
\item
\[(\ltn,\pi_{S}^{(A)})\sim\bigoplus_{J\in{\cal J}_{A,1}}P(J)
\oplus
\left\{
\bigoplus_{K\in{\cal J}_{A,\infty}}
P(K)\right\}^{\oplus\infty}\]
where ${\cal J}_{A,1}$ and ${\cal J}_{A,\infty}$
are in Definition \ref{defi:cycadef}.
%(ii)
\item
$\pi_{S}^{(A)}$ is multiplicity free if and only if
${\cal J}_{A,\infty}=\emptyset$.
Under this condition,
$\pi_{S}^{(A)}$ is irreducible if and only if
$\#{\cal J}_{A,1}=1$.
\end{enumerate}
\end{prop}
\pr
(i) By Lemma \ref{lem:abfsth},  the statement holds.

\noindent
(ii)
By (i), the first statement holds.
If $J=(j_{l})_{l=1}^{m}\in{\cal J}_{A,1}$,
then $j_{i}\ne j_{i^{'}}$ when $i\ne i^{'}$.
Therefore, any element in ${\cal J}_{A,1}$ is non periodic.
Hence the second statement holds.
\qedh

\noindent
If $A$ is full, that is, $\coa=\con$,
then $\pi_{S}^{(A)}$ is always $P(1)$ for each $N\geq 2$.

We show examples of $\pi_{S}^{(A)}$ by Proposition \ref{prop:stast}.
For this purpose, we use $s_{1},\ldots,s_{N}$ as 
canonical generators of $\coa$ and define operators
$t_{1},t_{2},k_{1},k_{2},k_{3}, u_{1},u_{2}$, $u_{3},u_{4}$ on $\ltn$ by
\[t_{i}e_{n}\equiv e_{2(n-1)+i},\quad k_{i}e_{n}\equiv e_{3(n-1)+i},
\quad u_{i}e_{n}\equiv e_{4(n-1)+i}\]
where $\{e_{n}:n\in {\bf N}\}$ is the canonical basis of $\ltn$.

Let
$A_{3}\equiv 
\matthr
{0}{1}{1}
{1}{0}{1}
{1}{1}{0}$ 
and $A_{4}\equiv 
\matthr{1}{0}{1}
{0}{1}{1}
{1}{1}{1}$.
The standard representation $(\ltn,\pi_{S}^{(A_{3})})$
of $\co{A_{3}}$ is given by
\[\pi^{(A_{3})}_{S}(s_{1})\equiv k_{1}(t_{1}k_{2}^{*}+t_{2}k_{3}^{*}),
\quad \pi^{(A_{3})}_{S}(s_{2})\equiv k_{2}(t_{1}k_{1}^{*}+t_{2}k_{3}^{*}),\]
\[\pi^{(A_{3})}_{S}(s_{3})\equiv k_{3}(t_{1}k_{1}^{*}+t_{2}k_{2}^{*}).\]
Then $\varphi_{A_{3}}$ in (\ref{eqn:varphia}) is given by
$\varphi_{A_{3}}:\{1,2,3\}\to \{1,2,3\};\quad
\varphi_{A_{3}}(1)=2$, $\varphi_{A_{3}}(2)=1$, $\varphi_{A_{3}}(3)=1$.
From this, $\varphi_{A_{3}}^{2}(1)=1$ and
${\cal J}_{A_{3}}={\cal J}_{A_{3},1}=\{(12)\}$.
Therefore \[(\ltn,\pi^{(A_{3})}_{S})\sim P(12).\]
%%%%%%%%%%%%%%%%%%%%%%%%%%%%
On the other hand, we have
\[
\pi^{(A_{4})}_{S}(s_{1})\equiv
k_{1}(t_{1}k_{1}^{*}+
t_{2}k_{3}^{*}),\quad
\pi^{(A_{4})}_{S}(s_{2})\equiv
k_{2}(t_{1}k_{2}^{*}+
t_{2}k_{3}^{*}),\quad
\pi^{(A_{4})}_{S}(s_{3})\equiv k_{3}.\]
Then ${\cal J}_{A_{4},1}=\{(1),(2)\}$,
$\pi^{(A_{4})}_{S}(s_{1})e_{1}=e_{1}$
and	$\pi^{(A_{4})}_{S}(s_{2})e_{2}=e_{2}$.
Therefore
\[(\ltn,\pi^{(A_{4})}_{S})\sim P(1)\oplus P(2).\]

Next, we show the decomposition of the standard representation
$\pi^{(A)}_{S}$ of $\coa$ for every element in $M_{2}(\{0,1\})$
without proof as follows:
{\tiny
\[
\begin{array}{c|c|c|c|c|c|c}
A&
\mattwo{1}{1}{1}{1},
\mattwo{1}{1}{1}{0}&
\mattwo{0}{1}{1}{1}&
\mattwo{1}{0}{1}{1}&
\mattwo{1}{1}{0}{1}&
\mattwo{1}{0}{0}{1}&
\mattwo{0}{1}{1}{0}\\
\hline
\pi^{(A)}_{S}& P(1)& P(12)& P(1)^{\oplus\infty}
& P(1)\oplus P(2)^{\oplus\infty}&
(P(1)\oplus P(2))^{\oplus\infty}& 
P(12)^{\oplus \infty}\\
\end{array}
\]
}
In this way,
the standard representation
of $\coa$ depends on $A$.
We show other examples as follows:
\[
\begin{array}{c|c}
A& \pi_{S}^{(A)}\\
\hline
\left(
\begin{array}{cccc}
0&1&0&1\\
0&1&0&1\\
1&1&0&1\\
0&1&1&1\\
\end{array}
\right)&
\begin{array}{l}
\pi_{S}^{(A)}(s_{1})=
u_{1}(t_{1}u_{2}^{*}+t_{2}u_{4}^{*}),\\
\pi_{S}^{(A)}(s_{2})=
u_{2}(t_{1}u_{2}^{*}+t_{2}u_{4}^{*}),\\
\pi_{S}^{(A)}(s_{3})=
u_{3}(k_{1}u_{1}^{*}+k_{2}u_{2}^{*}+
k_{3}u_{4}^{*}),\\
\pi_{S}^{(A)}(s_{4})=
u_{4}(k_{1}u_{2}^{*}+k_{2}u_{3}^{*}+
k_{3}u_{4}^{*}),\\
\pi_{S}^{(A)}\sim P(2).
\end{array}
\\
%%%%%%%%%%%%%%%%%%%%%%%%%%%%%
\hline
\left(
\begin{array}{cccc}
0&1&0&1\\
0&0&1&1\\
1&1&0&1\\
0&1&1&1\\
\end{array}
\right)&
\begin{array}{l}
\pi_{S}^{(A)}(s_{1})=
u_{1}(t_{1}u_{2}^{*}+t_{2}u_{4}^{*}),\\
\pi_{S}^{(A)}(s_{2})=
u_{2}(t_{1}u_{3}^{*}+t_{2}u_{4}^{*}),\\
\pi_{S}^{(A)}(s_{3})=
u_{3}(k_{1}u_{1}^{*}+k_{2}u_{2}^{*}+
k_{3}u_{4}^{*}),\\
\pi_{S}^{(A)}(s_{4})= u_{4}
(k_{1}u_{2}^{*}+k_{2}u_{3}^{*}+
k_{3}u_{4}^{*}),\\
\pi_{S}^{(A)}\sim P(123).
\end{array}
\\
\hline
\left(
\begin{array}{cccc}
0&0&1&1\\
1&0&1&1\\
0&1&0&1\\
0&1&1&1\\
\end{array}
\right)&
\begin{array}{l}
\pi_{S}^{(A)}(s_{1})=
u_{1}(t_{1}u_{3}^{*}+t_{2}u_{4}^{*}),\\
\pi_{S}^{(A)}(s_{2})=
u_{2}(k_{1}u_{1}^{*}+k_{2}u_{3}^{*}
+k_{3}u_{4}^{*}),\\
\pi_{S}^{(A)}(s_{3})=
u_{3}(t_{1}u_{2}^{*}+t_{2}u_{4}^{*}),\\
\pi_{S}^{(A)}(s_{4})=
u_{4}(k_{1}u_{2}^{*}+k_{2}u_{3}^{*}+
k_{3}u_{4}^{*}),\\
\pi_{S}^{(A)}\sim P(132).\\
\end{array}
\end{array}
\]

%%%%%%%%%%%%%%%%%%%%%%%%%%%%%%%%%%%%%%%%%%%%%%%
%
% subsection 7.3
%
\ssft{Shift representation}
\label{subsection:sevenththree}
For $A=(a_{ij})\in\mnz$, let $X_{A}\equiv \nset{\infty}_{A}$.
Define an $A$-\bfs\ $f=\{f_{i}\}_{i=1}^{N}$ on $X_{A}$ by
%
% Equation 8.1
%
\begin{equation}
\label{eqn:shift}
f_{i}:D(f_{i})\to R(f_{i});\quad 
f_{i}(j_{1},j_{2},\ldots)\equiv(i,j_{1},j_{2},\ldots),
\end{equation}
\[R(f_{i})\equiv\{(j_{1},j_{2},\ldots)\in X_{A}:j_{1}=i\},\quad
D(f_{i})\equiv \coprod_{j;a_{ij}=1}R(f_{j})\]
for $i\edot$.
The permutative representation $(l_{2}(X_{A}),\pi_{f})$ of $\coa$
by $f$ in (\ref{eqn:shift}) is called
the {\it shift representation} of $\coa$.
%
% Proposition 8.1
%
\begin{prop}
\label{prop:shift}
The following irreducible decomposition by $\pi_{f}$ holds:
\[l_{2}(X_{A})=\bigoplus_{J\in[1,\ldots,N]_{A}^{*}}{\cal H}_{J}
\quad\oplus\bigoplus_{K\in[1,\ldots,N]_{A}^{\infty}}{\cal K}_{K}\]
where
${\cal H}_{J}\equiv \overline{{\rm Lin}\langle\{e_{x}:x\in Y_{J^{\infty}}\}\rangle}$,
${\cal K}_{K}\equiv \overline{{\rm Lin}\langle\{e_{x}:x\in Y_{K}\}\rangle}$,
$\{e_{x}:x\in X_{A}\}$ is the canonical basis of $l_{2}(X_{A})$
and $Y_{L}\equiv \{L^{'}\in X_{A}:L\sim L^{'}\}$ for $L\in X_{A}$.
Furthermore
\[{\cal H}_{J}\sim P(J),\quad{\cal K}_{K}\sim P(K).\]
That is, any irreducible permutative representation
appears as a component of $(l_{2}(X_{A}),\pi_{f})$ once for all.
Especially, $(l_{2}(X_{A}),\pi_{f})$ is multiplicity free.
\end{prop}
%
% Proof
%
\pr
Let $X_{A,ep}$ be the set of all \evp\ elements in $X_{A}$
and $X_{A,nep}\equiv X_{A}\setminus X_{A,ep}$.
If $K\in X_{A,ep}$,	then there is $J\in[1,\ldots,N]_{A}^{*}$ 
such that $Y_{K}=Y_{J^{\infty}}$.
From this, $X_{A,ep}=\bigoplus_{J\in[1,\ldots,N]^{*}_{A}}Y_{J^{\infty}}$
and $f|_{Y_{J^{\infty}}}$ is $P(J)$. This implies the decomposition
$f|_{X_{A,ep}}\sim\bigoplus_{J\in[1,\ldots,N]^{*}_{A}}P(J)$.
Assume that $K=(k_{n})_{n\in{\bf N}}\in X_{A,nep}$ and 
$x_{n}\equiv (k_{n},k_{n+1},k_{n+2},\ldots)\in X_{A}$ for $n\geq 1$.
Then $X_{A,nep}=\bigoplus_{K\in[1,\ldots,N]^{\infty}_{A}}Y_{K}$,
$f|_{Y_{K}}$ is a cyclic $A$-\bfs\ on $Y_{K}$ and
$\{x_{n}\}_{n\in {\bf N}}$ is a chain of $f|_{Y_{K}}$ by $K$
in $Y_{K}$.
In consequence,
$f|_{X_{A,nep}}\sim\bigoplus_{K\in[1,\ldots,N]^{\infty}_{A}}P(K)$.
Because $X_{A}=X_{A,ep}\sqcup X_{A,nep}$, we have
the decomposition of $f$ into the direct sum of 
$\{P(J):J\in[1,\ldots,N]^{\#}_{A}\}$.
From this, the statement holds.
\qedh

\noindent
Proposition \ref{prop:shift} can be rewritten as follows:
\[(l_{2}(X_{A}),\pi_{f})\sim\bigoplus_{J\in[1,\ldots,N]_{A}^{\#}}P(J).\]

%%%%%%%%%%%%%%%%%%%%%%%%%%%%%%%%%%%%%%%%%%%%%%%	

%

\begin{thebibliography}{99}
%
\bibitem{AK1}Abe, M. and Kawamura, K.:
{\rm Recursive Fermion System in
Cuntz Algebra.\ I ---Embeddings of Fermion Algebra into Cuntz Algebra ---},
Comm.Math.Phys. {\bf 228}, 85-101 (2002),
%
\bibitem{AK3}Abe, M. and Kawamura, K.:
{\rm Pseudo Cuntz Algebra and Recursive FP Ghost System in String Theory},
Int. J. Mod. Phys. {\bf A18}, No. 4 (2003) 607-625.
% BBBBBB
%
\bibitem{BJ}Bratteli, O. and Jorgensen, P.E.T.:
{\rm Iterated function Systems and Permutation
Representations of the Cuntz algebra},
Memories Amer. Math. Soc. {\bf 139} (1999), no.663.
%
\bibitem{CK}Cuntz, J. and Krieger, W.:
{\rm A class of C$^{*}$-algebra and topological Markov chains},
Invent.Math., {\bf 56} (1980) 251-268.
%
%
\bibitem{DaPi2}Davidson, K.R. and Pitts, D.R.:
{\rm The algebraic structure of non-commutative analytic Toeplitz algebras},
Math.Ann. 311, 275-303 (1998),
{\rm Invariant subspaces and hyper-reflexivity for free semigroup algebras},
Proc. London Math. Soc. (3) 78 (1999) 401-430.
%
%%%%%%%%%%%%%%%%%%%%%%%%%%%
% Kawamura
%
%
\bibitem{PE01} Kawamura, K.:
{\rm Polynomial endomorphisms of the Cuntz algebras
arising from permutations. I ---General theory---},
Lett.Math.Phys., (2005) 71:149-158.
%
%
\bibitem{IWF01} Kawamura, K.:
{\rm Extensions of representations
of the CAR algebra to the Cuntz algebra $\co{2}$
---the Fock and the infinite wedge---},
J.Math.Phys., to appear.
%
%
\bibitem{PFO}Kawamura, K.:
{\rm The Perron-Frobenius operators, invariant measures
and representations of the Cuntz-Krieger algebras}, 
J.Math.Phys., to appear.
%
\end{thebibliography}
\end{document}